\documentclass[12pt]{article}
\usepackage{amssymb,amsmath,amsfonts}
\usepackage{relsize}
\usepackage[sort]{cite}
\usepackage{esvect}
\usepackage{txfonts}
\usepackage{graphicx}

\newtheorem{theorem}{Theorem}[section]
\newtheorem{definition}[theorem]{Definition}
\newtheorem{corollary}[theorem]{Corollary}
\newtheorem{proposition}[theorem]{Proposition}
\newtheorem{remark}[theorem]{Remark}
\newtheorem{lemma}[theorem]{Lemma}

\newtheorem{mlemma}[theorem]{Main Lemma}

\newtheorem{problem}[theorem]{Problem}
\newtheorem{question}[theorem]{Question}
\newtheorem{assumption}[theorem]{Assumption}

\newcommand {\Ec}      {{\mathcal E}}
\newcommand {\Fc}      {{\mathcal F}}

\newcommand {\Ic}      {{\mathcal I}}
\newcommand {\Jc}      {{\mathcal J}}
\newcommand {\Kc}      {{\mathcal K}}
\newcommand {\Lc}      {{\mathcal L}}
\newcommand {\Mc}      {{\mathcal M}}
\newcommand {\Nc}      {{\mathcal N}}

\newcommand {\Pc}      {{\mathcal P}}

\newcommand {\Sc}      {{\mathcal S}}

\newcommand {\Yc}      {{\mathcal Y}}
\newcommand {\Zc}      {{\mathcal Z}}
\newcommand {\R}       {{\mathbb R}}

\newcommand {\N}       {{\mathbb N}}
\newcommand {\Z}       {{\mathbb Z}}

\newcommand {\tH}      {\widetilde{H}}

\newcommand {\tS}      {\widetilde{S}}

\newcommand {\Jcw}     {\widetilde{\Jc}}
\newcommand {\Icw}     {\widetilde{\Ic}}

\newcommand {\tz}      {\tilde{z}}

\newcommand {\dw}      {\tilde{\delta}}

\newcommand {\RN}      {\R^n}
\newcommand {\ve}      {\varepsilon}
\newcommand {\LOPR}    {L^1_p(\R)}
\newcommand {\LMPR}    {L_{p}^{m}(\R)}
\newcommand {\WMP}     {W_{p}^{m}(\RN)}
\newcommand {\WMPR}    {W_{p}^{m}(\R)}
\newcommand {\LMIR}    {L^m_\infty(\R)}
\newcommand {\LMRN}    {L^m_\infty(\RN)}
\newcommand {\WMRN}    {W^m_\infty(\RN)}
\newcommand {\WTRT}    {W^2_\infty(\R^2)}

\newcommand {\LPR}     {L_p(\R)}
\newcommand {\CMR}     {C^{m}(\R)}

\newcommand {\intl}    {\int\limits}
\newcommand {\emp}     {\emptyset}
\newcommand {\PM}      {\Pc_{m}}
\newcommand {\PMO}     {\Pc_{m-1}}
\newcommand {\rl}      {r}
\newcommand {\VP}      {{\bf P}}
\newcommand {\brz}     {\bar{z}}
\newcommand {\SHF}     {\left(\Delta^mf\right)^\sharp_{E}}
\newcommand {\NMP}     {\Lc_{m,p}}
\newcommand {\TNMP}    {\widetilde{\Lc}_{m,p}}

\newcommand {\gmr}     {\gamma^{\sharp}(\LMIR)}
\newcommand {\gsho}    {\gamma^{\sharp}(L^1_\infty(\R))}
\newcommand {\gsht}    {\gamma^{\sharp}(L^2_\infty(\R))}

\newcommand {\NIN}     {\Lc}
\newcommand {\SH}      {{S\hspace*{-0.5mm}}}
\newcommand {\TSH}     {{\tS\hspace*{-0.4mm}}}

\newcommand {\LIM}     {\Zc}
\newcommand {\HH}      {\widehat{H}}

\newcommand {\xb}      {\bar{x}}

\newcommand {\WOP}     {\Fc_{m,E}^{(Wh)}}

\newcommand {\FOP}     {\Fc_{m,E}^{(Favard)}}
\newcommand {\KZ}      {K(z)}

\newcommand {\vkp}     {\varkappa}

\newcommand {\cm}      {\theta_m}
\newcommand {\CMM}     {\Theta_m}

\newcommand {\squ}     {\square\hspace{0.2mm}}
\newcommand {\smsk}    {\smallskip}
\newcommand {\msk}     {\medskip}

\newcommand {\mcup}    {\mathlarger{\cup}}
\newcommand {\mcap}    {\mathlarger{\cap}}
\newcommand {\pmed}    {\mathlarger{\prod}}

\newcommand {\smed}    {\mathlarger{\sum}}
\newcommand {\sbig}    {\mathlarger{\mathlarger{\sum}}}

\newcommand {\blbig}    {\mathlarger{\mathlarger{\bullet\hspace{0.2mm}}}}

\newcommand {\PME}     {\|\VP\|_{m,p,E}}
\newcommand {\VSH}     {\VP^\sharp_{m,E}}

\newcommand {\diam}    {\operatorname{diam}}
\newcommand {\dist}    {\operatorname{dist}}

\newcommand {\EXT}     {\operatorname{Ext}}

\newcommand {\VST}     {\vspace*{1mm}}
\newcommand {\bx}      {\hspace{10mm}$\Box$}
\newcommand {\rbx}     {\hspace{10mm}$\vartriangleleft$}

\newcommand {\BX}      {\hspace{10mm}\Box}
\newcommand {\nn}      {\nonumber}
\newcommand {\rf}[1]    {(\ref{#1})}      
\newcommand {\reff}[1] {\ref{#1}}         
\newcommand{\lbl}[1]      {\label{#1}}       
\newcommand{\be}          {\begin{eqnarray}}
\newcommand{\bel}[1]      {\begin{eqnarray} \label{#1}}
\newcommand{\ee}           {\end{eqnarray}}
\newcommand {\SECT}[2] {\section*{\centerline{\normalsize
{\bf #1}}} \setcounter{section}{#2}
\setcounter{theorem}{0}\setcounter{equation}{0}}
\begin{document}
\parindent 1em
\parskip 0mm
\medskip
\centerline{\large{\bf Extension criteria for homogeneous Sobolev spaces}}
\vspace*{5mm}
\centerline{\large{\bf of functions of one variable}} \vspace*{10mm}
\centerline{By~ {\sc Pavel Shvartsman}}
\vspace*{5 mm}
\centerline {\it Department of Mathematics, Technion - Israel Institute of Technology}\vspace*{2 mm}
\centerline{\it 32000 Haifa, Israel}\vspace*{2 mm}
\centerline{\it e-mail: pshv@technion.ac.il}
\vspace*{5 mm}
\renewcommand{\thefootnote}{ }
\footnotetext[1]{{\it\hspace{-6mm}Math Subject
Classification} 46E35\smallskip
\par This research was supported by Grant No 2014055 from the United States-Israel Binational Science Foundation (BSF).}

\begin{abstract} For each $p>1$ and each positive integer
$m$ we give intrinsic characterizations of the restriction of the homogeneous Sobolev space $\LMPR$
to an arbitrary closed subset $E$ of the real line.
\par We show that the classical one dimensional Whitney extension operator \cite{W2} is ``universal'' for the scale of $\LMPR$ spaces in the following sense: {\it for every $p\in(1,\infty]$ it provides almost optimal $L^m_p$-extensions of functions defined on $E$}. The operator norm of this extension operator is bounded by a constant depending only on $m$.
This enables us to prove several constructive $L^m_p$-extension criteria expressed in terms of $m^{\,\text{th}}$ order divided diffe\-rences of functions.
\end{abstract}
\renewcommand{\contentsname}{ }
\tableofcontents
\addtocontents{toc}{{\centerline{\sc{Contents}}}
\vspace*{5mm}\par}
\SECT{1. Introduction.}{1}
\addtocontents{toc}{~~~~1. Introduction.\hfill \thepage\par\VST}

\indent\par In this paper we characterize the restrictions of homogeneous Sobolev functions of one variable to an arbitrary closed subset of the real line. For each $m\in\N$ and each $p\in (1,\infty]$, we consider $\LMPR$,  the standard homogeneous Sobolev space on $\R$. We identify $\LMPR$ with the space of all real valued functions $F$ on $\R$ such that the $(m-1)$-th derivative $F^{(m-1)}$ is absolutely continuous on $\R$ and the weak $m$-th derivative $F^{(m)}\in L_p(\R)$. $\LMPR$ is seminormed by
$\|F\|_{\LMPR}= \|F^{(m)}\|_{L_p(\R)}$.
\smsk
\par  In this paper we study the following
\begin{problem}\lbl{PR-MAIN} {\em Let $p\in(1,\infty]$, $m\in\N$, and let $E$ be a closed subset of $\R$.
Let $f$ be a function on $E$. We ask two questions:\smallskip
\par {\it 1. How can we decide whether there exists a function $F\in \LMPR$ such that the restriction $F|_E$ of $F$ to $E$ coincides with $f$\,?}\smallskip
\par 2. Consider the $\LMPR$-seminorms of all functions $F\in\LMPR$ such that $F|_E=f$.  {\it How small can these seminorms be?}}
\end{problem}
\par We denote the infimum of all these seminorms by $\|f\|_{\LMPR|_E}$; thus
\bel{N-LMPR}
\|f\|_{\LMPR|_E}
=\inf \{\|F\|_{\LMPR}:F\in\LMPR, F|_{E}=f\}.
\ee
We refer to $\|f\|_{\LMPR|_E}$ as the {\it trace norm on $E$ of the function $f$} in $\LMPR$. This quantity provides the standard quotient space seminorm in {\it the trace space} $\LMPR|_{E}$ of all restrictions of $\LMPR$-functions to $E$, i.e., in the space 
$$
\LMPR|_{E}=\{f:E\to\R:\text{there exists}~~F\in
\LMPR\ \ \text{such that}\ \ F|_{E}=f\}.
$$
\par Whitney \cite{W2} completely solved an analog of part 1 of Problem \reff{PR-MAIN} for the space $C^m(\R)$. Whitney's extension construction \cite{W2} produces a certain {\it extension operator}
\bel{W-EXOP}
\WOP:C^m(\R)|_E\to C^m(\R)
\ee
which {\it linearly and continuously} maps the trace space  $C^m(\R)|_E$ into $C^m(\R)$. (See also Merrien \cite{Mer}.)
\par In fact the extension method developed by Whitney in \cite{W2} readily adapts to also provide a complete solution to Problem \reff{PR-MAIN} for the space $\LMIR$. Recall that $\LMIR$ can be identified with the space $C^{m-1,1}(\R)$ of all $C^{m-1}$-functions on $\R$ whose derivatives of order $m-1$ satisfy a Lipschitz condition. In particular, the method of proof and technique developed in \cite{W2} and \cite{Mer} lead us to the following well known description of the trace space $\LMIR|_E$: {\it A function $f\in\LMIR|_E$ if and only if the following quantity
$$
\NIN_{m,\infty}(f:E)=
\sup_{S\subset E,\,\,\# S=m+1}
|\Delta^mf[S]|
$$
is finite.} Here $\Delta^mf[S]$ denotes the {\it $m^{\,\text{th}}$ order divided difference} of $f$ on an $(m+1)$-point set $S$. 
\par Furthermore,
\bel{T-R1}
C_1\, \NIN_{m,\infty}(f:E)\le \|f\|_{\LMIR|_E}\le C_2\, \NIN_{m,\infty}(f:E)
\ee
where $C_1$ and $C_2$ are positive constants depending only on $m$. (Recall that $\Delta^mf[S]$ coincides with the coefficient of $x^m$ in the Lagrange polynomial of degree at most $m$ which agrees with $f$ on $S$. See Section 2.1 for other equivalent definitions of divided difference and their main properties.)
\smallskip
\par We refer the reader to \cite{Jon,SHE1} for further results in this direction.
\smallskip
\par There is an extensive literature devoted to a special case of Problem \reff{PR-MAIN} where $E$ consists of all the elements of a strictly increasing sequence $\{x_i\}_{i=\ell_1}^{\ell_2}$ (finite, one-sided infinite, or bi-infinite). We refer the reader to the papers of Favard \cite{Fav}, Chui, Smith, Ward \cite{ChS,CSW1,Sm-2}, de Boor \cite{deB2,deB4,deB5,deB6},
Fisher, Jerome \cite{FJ1}, Golomb \cite{Gol}, Jakimovski,  Russell \cite{JR}, Kunkle \cite{Kun2}, Pinkus \cite{P2}, Schoenberg \cite{Sch2} and references therein for numerous results in this direction and techniques for obtaining them.
\par In particular, for the space $\LMIR$ Favard \cite{Fav} developed a powerful linear extension method (very different from Whitney's method \cite{W2}) based on a certain delicate duality argument. Note that for any set $E$ as above and every $f:E\to \R$, Favard's extension operator $\FOP$ yields an extension of $f$ with {\it the smallest possible seminorm} in $\LMIR$. (Thus, $\|f\|_{\LMIR|_E}=\|\FOP(f)\|_{\LMIR}$ for every function $f$ defined on $E$.) Note also that  Favard's approach leads to the following slight refinement of \rf{T-R1}:
$$
\|f\|_{\LMIR|_E}\sim
\sup_{\ell_1\le i\le  \ell_2-m}
|\Delta^mf[x_i,...,x_{i+m}]|\,.
$$
See Section 7 for more details.
\smallskip
\par Modifying Favard's extension construction,  de Boor \cite{deB2,deB4,deB5} characterized the traces of $\LMPR$-functions to arbitrary sequences of points in $\R$. 
\begin{theorem}\lbl{DEBOOR}(\cite{deB4}) Let $p\in(1,\infty)$, and let $\ell_1,\ell_2\in\Z\cup\{\pm\infty\}$, $\ell_1+m\le\ell_2$. Let $f$ be a function defined on a strictly increasing sequence of points $E=\{x_i\}_{i=\ell_1}^{\ell_2}$. Then
$f\in\LMPR|_E$ if an only if the following quantity
$$
\TNMP(f:E)=\,\left(\smed_{i=\ell_1}^{\ell_2-m}\,\,
(x_{i+m}-x_i)\,|\Delta^mf[x_i,...,x_{i+m}]|^p
\right)^{\frac1p}
$$
is finite. Furthermore,
$\|f\|_{\LMPR|_E}\sim \TNMP(f:E)$ with constants of equivalence depending only on $m$.
\end{theorem}
\par For a special case of this result, for sequences satisfying some global mesh ratio restrictions, see Golomb \cite{Gol}. See also Est\'evez \cite{Es} for an alternative proof of Theorem \reff{DEBOOR} for $m=2$.
\par Using a certain limiting argument, Golomb \cite[Theorem 2.1]{Gol} showed that Problem \reff{PR-MAIN} for $\LMPR$ and an {\it arbitrary} set $E\subset\R$ can be reduced to the same problem, but for arbitrary {\it finite} sets $E$. More specifically, his result (in an equivalent form) provides the following formula for the trace norm in $\LMPR|_{E}$:
$$
\|f\|_{\LMPR|_E}=\sup\{\,\|f|_{E'}\|_{\LMPR|_{E'}}: E'\subset E,\# E'<\infty\}.
$$
\par Let us remark that, by combining this formula with de Boor's Theorem \reff{DEBOOR}, we can obtain the following description of the trace space $\LMPR|_E$ for an {\it arbitrary} closed set $E\subset\R$.
\begin{theorem}\lbl{MAIN-TH}(Variational extension criterion for  $\LMPR$-traces)~ Let $p\in(1,\infty)$ and let $m$ be a positive integer. Let $E\subset\R$ be a closed set containing at least $m+1$ points. A function $f:E\to\R$ can be extended to a function $F\in\LMPR$ if and only if the following quantity
\bel{NR-TR}
\NMP(f:E)=\,\sup_{\{x_0,...,x_n\}\subset E}
\,\,\,\left(\smed_{i=0}^{n-m}
\,\,
(x_{i+m}-x_i)\,|\Delta^mf[x_i,...,x_{i+m}]|^p
\right)^{\frac1p}
\ee
is finite. Here the supremum is taken over all all integers $n\ge m$ and all strictly increasing sequences $\{x_0,...,x_n\}\subset E$ of $n$ elements. Furthermore,
\bel{NM-CLC}
\|f\|_{\LMPR|_E}\sim \NMP(f:E).
\ee
The constants of equivalence \rf{NM-CLC} depend only on $m$.
\end{theorem}
\par In the present paper we give a {\it direct and explicit proof} of Theorem \reff{MAIN-TH} which {\it does not use any limiting argument}. Actually we show, perhaps surprisingly, that {\it the very same Whitney extension ope\-rator $\WOP$} (see \rf{W-EXOP}) which was introduced in \cite{W2} for characterization of the trace space $C^{m}(\R)|_E$, {\it provides almost optimal extensions of functions belonging to}  $\LMPR|_E$ {\it for every $p\in(1,\infty]$}.
\smallskip
\par We also give another characterization of the trace space $\LMPR|_E$ expressed in terms of $L^p$-norms of certain kinds of ``sharp maximal functions'' which are defined as follows:
\par For each $m\in\N$, each closed set $E\subset\R$ with $\# E>m$, and each function $f:E\to\R$ we let $\SHF$ denote the maximal function associated with $f$ which is given by
\bel{SH-F}
\SHF(x)=
\,\sup_{\substack {\{x_0,...,x_m\}
\subset E\smallskip\\ x_0<x_1<...<x_m}}\,\,\,
\frac{|\,\Delta^{m-1}f[x_0,...,x_{m-1}]-
\Delta^{m-1}f[x_1,...,x_{m}]|}{|x-x_0|+|x-x_m|},~~~x\in\R\,.
\ee
\begin{theorem} \lbl{R1-CR} Let $p\in (1,\infty)$, $m\in\N$, and let $f$ be a function defined on a closed set $E\subset\R$. The function $f\in \LMPR|_E$ if and only if $\SHF\in \LPR$. Furthermore,
$$
\|f\|_{\LMPR|_E}\sim \|\SHF\|_{\LPR}
$$
with the constants in this equivalence depending only on $m$ and $p$.
\end{theorem}
\par Note that
\bel{SH-F-EQ}
\SHF(x)\,\le \,\sup_{S\subset E,\,\,\# S=m+1}\,\,
\frac{|\Delta^mf[S]|\diam S}
{\diam (\{x\}\cup S)}\le 2\,\SHF(x)~~~\text{for all}~~~x\in\R\,.
\ee
(See property \rf{D-IND} below.) This inequality, together with Theorem \reff{R1-CR} and the definition in \rf{SH-F} now imply the following explicit formulae for the trace seminorm of a function in the space $\LMPR|_E$: 
\be
\|f\|_{\LMPR|_E}&\sim&
\left\{\intl\limits_{\R}
\,\sup_{\substack {\{x_0,...,x_m\}
\subset E\smallskip\\ x_0<x_1<...<x_m}}
\,\,\frac{|\,\Delta^{m-1}f[x_0,...,x_{m-1}]-
\Delta^{m-1}f[x_1,...,x_{m}]|^{\,p}}{|x-x_0|^p+|x-x_m|^{p}}
\,dx\right\}^{\frac{1}{p}}
\nn\\
&\sim&
\left\{\,\intl\limits_{\R}\sup_{S\subset E,\,\,\# S=m+1}\,\left(\frac{|\Delta^mf[S]|\diam S}{\diam (\{x\}\cup S)}\right)^pdx\right\}^{\frac{1}{p}}.\nn
\ee
\par We feel a strong debt to the remarkable papers of Calder\'{o}n and Scott \cite{C1,CS} which are devoted to characterization of Sobolev spaces on $\RN$ in terms of  classical sharp maximal functions. These papers motivated us to formulate and subsequently prove Theorem \reff{R1-CR}.
\par For analogs of Theorems \reff{MAIN-TH} and \reff{R1-CR} for the space $L^1_p(\RN)$, $n\in\N$, $n<p<\infty$, we refer the reader to \cite{Sh2,Sh5}.
\medskip
\par Our next new result, Theorem \reff{LIN-OP} below, states that there exists a solution to Problem \reff{PR-MAIN} which depends {\it linearly} on the initial data, i.e., the functions defined on $E$.
\begin{theorem} \lbl{LIN-OP} For every closed subset $E\subset\R$, every $p>1$  and every $m\in\N$ there exists a conti\-nuous linear extension operator which maps the trace space $\LMPR|_E$ into $\LMPR$. Its operator norm is
bounded by a constant depending only on $m$.
\end{theorem}
\par Let us recall something of the history of the previous results which led us to Theorem \reff{LIN-OP}. We know that for each closed $E\subset\R$ the Whitney extension operator $\WOP$ \cite{W2} maps $\LMIR|_E$ into $\LMIR$ with the operator norm $\|\WOP\|$ bounded by a constant depending only on $m$. As we have mentioned above, if $E$ is a sequence of points in $\R$, Favard's linear extension operator also maps $\LMIR|_E$ into $\LMIR$, but with the operator norm  $\|\FOP\|=1$.
\par For $p\in(1,\infty)$ and an arbitrary sequence $E\subset\R$ Theorem \reff{LIN-OP} follows from \cite[Section 4]{deB4}. Luli \cite{L} gave an alternative proof of Theorem \reff{LIN-OP} for the space $\LMPR$ and a finite set $E$. In the multidimensional case the existence of corresponding linear continuous extension operators for the Sobolev spaces $L^m_p(\RN)$, $n<p<\infty$, was proven in \cite{Sh2} ($m=1$, $n\in\N$, $E\subset\RN$ is arbitrary), \cite{Is} and \cite{Sh3} ($m=2, n=2$, $E\subset\R^2$ is finite), and \cite{FIL} (arbitrary $m,n\in\N$ and an arbitrary $E\subset\RN$). For the case $p=\infty$ see \cite{BS2} ($m=2$) and \cite{F3,F-IO} ($m\in\N$).
\medskip
\par In a forthcoming paper \cite{Sh-2018} we will present a solution to an analog of Problem \reff{PR-MAIN} for the normed Sobolev space $\WMPR$.
\medskip
\par Let us briefly describe the structure of the present paper and the main ideas of our approach.
\smallskip
\par First we note that the equivalence \rf{NM-CLC} is not trivial even in the simplest case, i.e., for $E=\R$; in this case \rf{NM-CLC} tells us that for every $f\in\LMPR$ and every $p\in(1,\infty]$
\bel{RNM-LMP}
\|f\|_{\LMPR}\sim \NMP(f:\R)
\ee
with constants depending only on $m$. In other words, the quantity $\NMP(\cdot:\R)$ provides an equivalent seminorm on $\LMPR$. This characterization of the space $\LMPR$ is known in the literature; see F.~Riesz \cite{R} ($m=1$ and $1<p<\infty$), Schoenberg \cite{Sch2} ($p = 2$ and $m\in\N$), and Jerome and Schumaker \cite{JS} (arbitrary $m\in\N$ and $p\in(1,\infty)$). Of course, the equivalence \rf{RNM-LMP} implies {\it the necessity part} of Theorem \reff{MAIN-TH}. Nevertheless, for the reader's convenience, in Section 2.2 we give a short direct proof of this result (together with the proof of the necessity part of Theorem \reff{R1-CR}).
\smsk
\par In Section 3 we recall the Whitney extension method \cite{W2} for functions of one variable. We prove a series of auxiliary statements which enable us to adapt Whitney's construction to extension of $\LMPR$-functions. We then use this extension technique and a criterion for extension of Sobolev jets \cite{Sh5} to help us prove the sufficiency part of Theorem \reff{R1-CR}. (See Section 3.4.)
\smallskip
\par The sufficiency part of Theorem \reff{MAIN-TH} is proven in Sections 4-6. One of the main ingredients of this proof is Theorem \reff{JET-V}, a refinement of the extension criterion given in Theorem \reff{JET-S}. Another important ingredient of the proof of the sufficiency is Main Lemma \reff{X-SXN}  which provides a certain controlled transition from Hermite polynomials of a function to its Lagrange polynomials. See Section 5.
\par In Section 6, with the help of these results, Theorem \reff{JET-V} and Main Lemma \reff{X-SXN}, we complete the proof of the sufficiency part of Theorem \reff{MAIN-TH}.
\msk
\par In Section 7 we discuss the dependence on $m$ of the constants $C_1,C_2$ in inequality \rf{T-R1}. We interpret this inequality as a particular case of {\it the Finiteness Principle for traces of smooth functions}. (See Theorem \reff{FP-LM}). We refer the reader to \cite{BS1,BS3,F2,F-Bl,Sh-2008} and references therein for numerous results related to the Finiteness Principle.
\par For the space $\LMIR$ the Finiteness Principle is equivalent to the following statement: there exists a constant $\gamma=\gamma(m)$ such that for every closed set $E\subset\R$ and every $f\in\LMIR|_E$ the following inequality
\bel {FN-11}
\|f\|_{\LMIR|_E}\le \gamma
\,\sup_{S\subset E,\,\,\# S= m+1}\,\|f|_S\|_{\LMIR|_S}
\ee
holds. We can express this result by stating that the number $m+1$ is a {\it finiteness number} for the space $\LMIR$. We also refer to any constant $\gamma$
which satisfies \rf{FN-11} as {\it a multiplicative finiteness constant} for the space $\LMIR$. In this context we let $\gmr$ denote the infimum of all multiplicative finiteness constants for $\LMIR$ for the finiteness number $m+1$.
\par One can easily see that $\gsho=1$. In Theorem \reff{G-SH} we show that
\bel {FC-IN}
\gsht=2~~~\text{and}~~~
(\pi/2)^{m-1}<\gmr<(m-1)\,9^m~~~\text{for every}~~~m>2\,.
\ee
\par The proof of \rf{FC-IN} relies on results of Favard \cite{Fav} and de Boor \cite {deB4,deB5} devoted to calculation of certain extension constants for the space $\LMIR$. See Section 7 for more details.
\smsk
\par Readers might find it helpful to also
consult a much more detailed version of this paper posted on the arXiv \cite{Sh-LV-2018}.

\bigskip
\par {\bf Acknowledgements.} I am very thankful to M. Cwikel for useful suggestions and remarks.
\par I am grateful to Charles Fefferman, Bo'az Klartag and Yuri Brudnyi for valuable conversations. The results of this paper were presented at the 11th Whitney Problems Workshop, Trinity College Dublin, Dublin, Ireland. I am very thankful to all participants of this conference for stimulating discussions and valuable advice.
\SECT{2. Main Theorems: necessity.}{2}
\addtocontents{toc}{2. Main Theorems: necessity.\hfill \thepage\par\VST}

\indent\par 

\par Let us fix some notation. Throughout the paper $C,C_1,C_2,...$ will be generic posi\-tive constants which depend only on $m$ and $p$. These symbols may denote different constants in different occurrences. The dependence of a constant on certain parameters is expressed by the notation $C=C(m)$, $C=C(p)$ or $C=C(m,p)$.
Given constants $\alpha,\beta\ge 0$, we write $\alpha\sim \beta$ if there is a constant $C\ge 1$ such that $\alpha/C\le \beta\le C\,\alpha$.
\par Given a measurable set $A\subset \R$, we let  $\left|A\right|$ denote the Lebesgue measure of $A$.
If $A\subset\R$ is finite, by $\#A$ we denote the number of elements of $A$.
\par Given $A,B\subset \R$, let
$$
\diam A=\sup\{\,|\,a-a'\,|:~a,a'\in A\}~~~\text{and}~~~ \dist(A,B)=\inf\{\,|\,a-b\,|:~a\in A, b\in B\}.
$$
For $x\in \R$ we also set $\dist(x,A)=\dist(\{x\},A)$. The notation
$$
A\to x~~~\text{will mean that}~~~
\diam(A\cup \{x\})\to 0.
$$
\par Given $M>0$ and a family $\Ic$  of intervals in $\R$ we say that {\it covering multiplicity} of $\Ic$ is bounded by $M$ if every point $x\in\R$ is covered by at most $M$ intervals from $\Ic$.
\par Given a function $g\in L_{1,loc}(\R)$ we let $\Mc[g]$ denote the Hardy-Littlewood maximal function of $g$:
\bel{HL-M}
\Mc[g](x)=\sup_{I\ni x}\frac{1}{|I|}\intl_I|g(y)|dy,~~~~x\in\R.
\ee
Here the supremum is taken over all closed intervals $I$ in $\R$ containing $x$.
\par By $\PM$ we denote the space of all polynomials of degree at most $m$ defined on $\R$. Finally, given a nonnegative integer $k$, a $(k+1)$-point set $S\subset\R$ and a function $f$ on $S$, we let $L_S[f]$ denote the Lagrange polynomial of degree at most $k$ interpolating $f$ on $S$; thus
$$
L_S[f]\in \Pc_k~~~~\text{and}~~~~L_S[f](x)=f(x)~~~
\text{for every}~~~x\in S.
$$

\bigskip\bigskip
\par {\bf 2.1. Divided differences: main properties.}
\medskip
\addtocontents{toc}{~~~~2.1. Divided differences: main properties. \hfill \thepage\par}
\par In this section we recall several useful properties of the divided differences of functions. We refer the reader to  \cite[Ch. 4, \S 7]{DL}, \cite[Section 1.3]{FK-88} and \cite[Section 2.1]{Sh-LV-2018} for the proofs of these properties.\smallskip
\par Everywhere in this section $k$ is a nonnegative integer and $S=\{x_0,...,x_k\}$ is a $(k+1)$-point subset of $\R$. In $(\bigstar 1)$-$(\bigstar 3)$ by $f$ we denote a function defined on $S$.
\msk
\par Then the following properties hold:
\smallskip
\par $(\bigstar 1)$~ $\Delta^0f[S]=f(x_0)$ provided $S=\{x_0\}$ is a singleton.
\medskip
\par $(\bigstar 2)$~ If $k\in\N$ then
\bel{D-IND}
\Delta^kf[S]=\Delta^{k}f[x_0,x_1,...,x_{k}]
=\left(\Delta^{k-1}f[x_1,...,x_{k}]
-\Delta^{k-1}f[x_0,...,x_{k-1}]\right)/(x_k-x_0).
\ee
\par Furthermore,
$$
\Delta^kf[S]=
\smed_{i=0}^k\,\,\frac{f(x_i)}{\omega'(x_i)}=
\smed_{i=0}^k\,\,\frac{f(x_i)}
{\prod\limits_{j\in\{0,...,k\},j\ne i}(x_i-x_j)}
$$
where $\omega(x)=(x-x_0)...(x-x_k)$.
\medskip
\par $(\bigstar 3)$~ We recall that $L_S[f]$ denotes the Lagrange polynomial of degree at most $k=\#S-1$ interpolating $f$ on $S$. Then the following equality
\bel{D-LAG}
\Delta^{k}f[S]=\frac{1}{k!}\,L^{(k)}_S[f]
\ee
holds. Thus,
$\Delta^{k}f[S]=A_k~~~\text{where}~~A_k~~\text{is the coefficient of}~~x^k~~\text{of the polynomial}~~L_S[f]$.
\msk
\par $(\bigstar 4)$~ Let $k\in\N$, and let $x_0=\min\{x_i:i=0,...,k\}$ and $x_k=\max\{x_i:i=0,...,k\}$. Then for every function $F\in C^k[x_0,x_k]$ there exists $\xi\in [x_0,x_k]$ such that
\bel{D-KSI}
\Delta^{k}F[x_0,x_1,...,x_{k}]
=\frac{1}{k!}\,F^{(k)}(\xi)\,.
\ee

\par $(\bigstar 5)$~ Let $x_0<x_1<...<x_k$, and let $F$ be a function on $[x_0,x_k]$ with absolutely continuous derivative of order $k-1$. Then
\bel{DVD-IN}
|\Delta^{k}F[S]|\le\frac{1}{(k-1)!}\cdot\frac{1}{x_k-x_0}
\,\intl_{x_0}^{x_k}\,
|F^{(k)}(t)|\,dt.
\ee
\par Furthermore, for every $\{x_0,...,x_m\}\subset\R$, $x_0<...<x_m$, and every $F\in\LMIR$ the following inequality
\bel{F-LMIR}
m!\,|\Delta^mF[x_0,...,x_m]|\le\|F\|_{\LMIR}
\ee
holds.


\bigskip\bigskip
\par {\bf 2.2. Proofs of the necessity parts of the main theorems.}
\medskip
\addtocontents{toc}{~~~~2.2. Proofs of the necessity parts of the main theorems.\hfill \thepage\VST\par}
\par {\it (Theorem \reff{MAIN-TH}: Necessity)} Let $1<p<\infty$ and let $f\in\LMPR|_E$. Let $F\in\LMPR$ be an arbitrary function such that $F|_E=f$. Let $n\ge m$ and let $\{x_0,...,x_n\}\subset E$, $x_0<...<x_n$. From \rf{DVD-IN}, for every $i, 0\le i\le n-m$, we have
\be
A_i&=&(x_{i+m}-x_i)\,|\Delta^m f[x_i,...,x_{i+m}]|^p
=
(x_{i+m}-x_i)\,|\Delta^mF[x_i,...,x_{i+m}]|^p\nn\\
&\le&
(x_{i+m}-x_i)\cdot\frac{1}{((m-1)!)^p}\cdot
\,\frac{1}{x_{i+m}-x_i}
\,\intl_{x_i}^{x_{i+m}}\,|F^{(m)}(t)|^p\,dt=
\frac{1}{((m-1)!)^p}
\,\intl_{x_i}^{x_{i+m}}\,|F^{(m)}(t)|^p\,dt\,.\nn
\ee
Hence,
$$
A=\smed_{i=0}^{n-m}\,A_i
\le
\frac{1}{((m-1)!)^p}\,\smed_{i=0}^{n-m}\,
\,\intl_{x_i}^{x_{i+m}}\,|F^{(m)}(t)|^p\,dt\,.
$$
\par Clearly, the covering multiplicity of the family $\{(x_i,x_{i+m}):i=0,...,n-m\}$ of open intervals is bounded by $2m$, so that
$$
A\le
\,\frac{2m}{((m-1)!)^p}
\,\intl_{x_0}^{x_{n}}\,|F^{(m)}(t)|^p\,dt\le
\frac{2m}{((m-1)!)^p}\,\|F\|_{\LMPR}^p\le
4^p\,\|F\|_{\LMPR}^p\,.
$$
\par This inequality together with definition \rf{NR-TR} implies that $\NMP(f:E)\le 2\,\|F\|_{\LMPR}$.
\par Finally, taking the infimum in the right hand side of this inequality over all functions $F\in\LMPR$ such that $F|_E=f$, we conclude that
$$
\NMP(f:E)\le 2\,\|f\|_{\LMPR|_E}
$$
proving the necessity part of Theorem \reff{MAIN-TH}.\bx
\bigskip
\par {\it (Theorem \reff{R1-CR}: Necessity)} Let $p\in(1,\infty)$ and let $f\in\LMPR|_E$. Let $F\in\LMPR$ be an arbitrary function such that $F|_E=f$. Let $S=\{x_0,...,x_m\}\subset E$, $x_0<...<x_m$, and let $x\in\R$.
\par From \rf{D-IND} and \rf{DVD-IN}, we have
\be
B&=&
\frac{|\,\Delta^{m-1}f[x_0,...,x_{m-1}]-
\Delta^{m-1}f[x_1,...,x_{m}]|}{|x-x_0|+|x-x_m|}
=
\frac{|\,\Delta^{m-1}F[x_0,...,x_{m-1}]-
\Delta^{m-1}F[x_1,...,x_{m}]|}{|x-x_0|+|x-x_m|}\nn\\
&=&
\frac{|\Delta^{m}F[x_0,...,x_{m}]|\,(x_m-x_0)}
{|x-x_0|+|x-x_m|}
\le
\frac{1}{(m-1)!\,(|x-x_0|+|x-x_m|)}
\intl_{x_0}^{x_{m}}\,|F^{(m)}(t)|\,dt\,.\nn
\ee
\par Let $I$ be the smallest closed interval containing $S$ and $x$. Clearly, $|I|\le |x-x_0|+|x-x_m|$ and $I\supset[x_0,x_m]$. Hence,
$$
B\le
\frac{1}{(m-1)!}\,\frac{1}{|I|}
\intl_I\,|F^{(m)}(t)|\,dt\le \,\Mc[F^{(m)}](x)\,.
~~~~\text{See \rf{HL-M}.}
$$
\par This inequality together with definition \rf{SH-F} implies that $\SHF(x)\le\,\Mc[F^{(m)}](x)$ on $\R$. Hence,
$$
\|\SHF\|_{\LPR}
\le\,\|\Mc[F^{(m)}]\|_{\LPR},
$$
so that, by the Hardy-Littlewood maximal theorem,
$$
\|\SHF\|_{\LPR}
\le\,C(p)\,\|F^{(m)}\|_{\LPR}=\,C(p)\,\|F\|_{\LMPR}\,.
$$
\par Taking the infimum in the right hand side of this inequality over all functions $F\in\LMPR$ such that $F|_E=f$, we finally obtain the required inequality
$$
\|\SHF\|_{\LPR}\le\,C(p)\,\|f\|_{\LMPR|_E}\,.
$$
\par The proof of the necessity part of Theorem \reff{R1-CR} is complete.\bx
\bigskip

\SECT{3. The Whitney extension method in $\R$ and traces of Sobolev functions.}{3}
\addtocontents{toc}{3. The Whitney extension method in $\R$ and traces of Sobolev functions.\hfill \thepage\par\VST}
\indent\par 
\par In this section we prove the sufficiency part of Theorem \reff{R1-CR}.
\par Given a function $F\in \CMR$ and $x\in\R$, we let 
$$
T^m_x[F](y)=\smed_{k=0}^m\,\,\frac{1}{k!}\, F^{(k)}(x)(y-x)^{k},~~~~y\in\R,
$$
denote the Taylor polynomial of $F$ of degree $m$ at $x$.
\smsk
\par Let $E$ be a closed subset of $\R$, and let $\VP=\{P_x: x\in E\}$ be a family of polynomials of degree at most $m$ indexed by points of $E$. (Thus $P_x\in \PM$ for every $x\in E$.) Following \cite{F9}, we refer to $\VP$ as {\it a Whitney $m$-field defined on $E$}.
\par We say that a function $F\in \CMR$ {\it agrees with the Whitney $m$-field $\VP=\{P_x: x\in E\}$ on $E$}, if $T^{m}_x[F]=P_x$ for each $x\in E$. In that case we also refer to $\VP$ as the Whitney $m$-field on $E$ {\it generated by $F$} or as {\it the $m$-jet generated by $F$.} We define the $L^m_p$-``norm'' of the $m$-jet $\VP=\{P_x: x\in E\}$ by
\bel{N-VP}
\PME=\inf\left\{\|F\|_{\LMPR}:F\in \LMPR,\, T^{m-1}_x[F]=P_x~~\text{for every}~~x\in E\right\}.
\ee
\par We prove the sufficiency part of Theorem \reff{R1-CR} in two steps. At the first step, given $m\in\N$ we construct a linear operator which to every function $f$ on $E$ assigns a certain  Whitney $(m-1)$-field
\bel{PME-1}
\VP^{(m,E)}[f]=\{P_x\in\PMO:x\in E\}
\ee
such that $P_x(x)=f(x)$ for all $x\in E$. We produce  $\VP^{(m,E)}[f]$ by a slight modification of Whitney's extension construction \cite{W2}. See also \cite{Mer,FK-88,FK-93,KM} where similar constructions have been used for characterization of traces of $\LMIR$-functions.
\smallskip
\par At the second step of the proof we show that for every $p\in(1,\infty)$ and every function $f:E\to\R$ such that $\SHF\in\LPR$ (see \rf{SH-F}) the following inequality
\bel{TR-PW}
\|\VP^{(m,E)}[f]\|_{m,p,E}\le C(m,p)\|\SHF\|_{\LPR}
\ee
holds. One of the main ingredients of the proof of \rf{TR-PW} is a trace criterion for jets generated by Sobolev functions. See Theorem \reff{JET-S} below.
\bigskip\medskip
\par {\bf 3.1. Interpolation knots and their properties.}
\medskip
\addtocontents{toc}{~~~~3.1. Interpolation knots and their properties. \hfill \thepage\par}
\par Let $E\subset\R$ be a closed subset, and let $k$ be a non-negative integer, $k\le\#E$. Following \cite{W2} (see also \cite{KM,Mer}), given $x\in E$ we construct an important ingredients of our extension procedure, a finite set $Y_k(x)\subset E$, which, in a certain sense, is ``well concentrated'' around $x$. This set provides interpolation knots for Lagrange and Hermite polynomials which we use in our modification of the Whitney extension method.
\par We will need the following notion. Let $A$ be a nonempty {\it finite subset} of $E$, $A\ne E$. Suppose that  $A$ contains {\it at most one limit point of $E$}. We assign to $A$ a point $a_E(A)\in E$ in the {\it closure} of $E\setminus A$ {\it having the minimal distance} to $A$. More specifically:
\smallskip
\par (i) If $A$ does not contain limit points of $E$, the set $E\setminus A$ is non-empty and closed, so that in this case $a_E(A)$ is a point nearest to $A$ on $E\setminus A$. Clearly, in this case  $a_E(A)\notin A$;
\smallskip
\par (ii) Suppose there exists a (unique) point $a\in A$ which is a limit point of $E$. In this case we set $a_E(A)=a$.
\medskip
\par Note that in both cases
$$
\dist(a_E(A),A)=\dist(A,E\setminus A).
$$
\par Now, let us construct a family of points $\{y_0(x),y_1(x),...,y_{n_k(x)}\}$ in $E$, $0\le n_k(x)\le k$, using the following inductive procedure.
\par First, we put $y_0(x)=x$ and $Y_0(x)=\{y_0(x)\}$. If $k=0$, we put $n_k(x)=0$, and stop.
\par Suppose that $k>0$. If $y_0(x)=x$ is a {\it limit point} of $E$, we again put $n_k(x)=0$, and stop. If $y_0(x)$ is an {\it isolated point of $E$}, we continue the procedure.
\par We define a point $y_1(x)\in E$ by
$y_1(x)=a_E(Y_0(x))$, and set $Y_1(x)=\{y_0(x),y_1(x)\}$. If $k=1$ or $y_1(x)$ is a limit point of $E$, we put $n_k(x)=1$, and stop.
\par Let $k>1$ and $y_1(x)$ is an {\it isolated point of $E$}. In this case we put
$$
y_2(x)=a_E(Y_1(x))~~~\text{and}~~~ Y_2(x)=\{y_0(x),y_1(x),y_2(x)\}.
$$
\par If $k=2$ or $y_2(x)$ is a limit point of $E$, we set
$n_k(x)=2$, and stop. But if $k>2$ and $y_2(x)$ is an isolated point of $E$, we continue the procedure and  define $y_3$, etc.
\par At the $j$-th step of this algorithm we obtain a $j+1$-point set $Y_j(x)=\{y_0(x),...,y_j(x)\}$.
\smsk
\par If $j=k$ or $y_j(x)$ is a limit point of $E$, we put $n_k(x)=j$ and stop. But if $j<k$ and $y_j(x)$ is an isolated point of $E$, we define a point $y_{j+1}(x)$ and a set $Y_{j+1}(x)$ by the formulae
\bel{DF-JY1}
y_{j+1}(x)=a_E(Y_{j}(x))~~~\text{and}~~~
Y_{j+1}(x)=\{y_0(x),...,y_{j}(x),y_{j+1}(x)\}\,.
\ee
\par Clearly, for a certain
$n=n_k(x)$,~ $0\le n\le k$, the procedure stops.
This means that either $n=k$ or, when $n<k$, the points $y_0(x),...,y_{n-1}(x)$ are {\it isolated points} of $E$, but
\bel{YN}
y_n(x)~~~\text{is a limit point of}~~~E\,.
\ee
\par We also introduce points $y_j(x)$ and sets $Y_j(x)$ for $n_k(x)\le j\le k$ by letting
\bel{YJ-NX}
y_j(x)=y_{n_k(x)}(x)~~~\text{and}~~~
Y_j(x)=Y_{n_k(x)}(x).
\ee
\par Note that, given $x\in E$ the definitions of points $y_j(x)$ and the sets $Y_j(x)$ {\it do not depend on $k$}, i.e, $y_j(x)$ is the same point and $Y_j(x)$ is the same set for every  $k\ge j$. This is immediate from \rf{YJ-NX}.
\par In the next three lemmas, we describe several important properties of the points $y_j(x)$ and sets $Y_j(x)$.
\begin{lemma}\lbl{L-YP} Given $x\in E$, the points $y_j(x)$ and the sets $Y_j(x)$, $0\le j\le k$, have the following properties:\medskip 
\par (a). $y_0(x)=x$ and $y_j(x)=a_E(Y_{j-1}(x))$ for every $1\le j\le k$;\smallskip
\par (b). Let $n=n_k(x)\ge 1$. Then  $y_0(x),...,y_{n-1}(x)$ are isolated points of $E$.
\smsk
\par Furthermore, if $y\in Y_n(x)$ and $y$ is a limit point of $E$, then $y=y_n(x)$.
\par In addition, if~ $0<n<k$, then $y_n(x)$ is a limit point of $E$;
\smallskip
\par (c). $\#Y_j(x)=\min\{j,n_k(x)\}+1$ for every $0\le j\le k$;
\smallskip
\par (d). For every $j=0,...,k,$ we have
\bel{MNM-Y}
[\min Y_j(x),\max Y_j(x)]\cap E= Y_j(x).
\ee
Furthermore, the point $y_j(x)$ is either minimal or maximal point of the set $Y_j(x)$.
\end{lemma}
\par {\it Proof.} Properties {\it (b)-(d)} are immediate from the definitions of the points $y_j(x)$ and sets $Y_j(x)$.
\par Let us prove {\it (a)}. We know that $y_0(x)=x$ and, thanks to \rf{DF-JY1},  $y_j(x)=a_E(Y_{j-1}(x))$ for every $j=1,...,n_k(x)$. If $n_k(x)<j\le k$, then, by \rf{YJ-NX},
$Y_j(x)=Y_{n_k(x)}(x)$ for every $j,~n_k(x)\le j\le k$.
\smsk
\par On the other hand, since $n_k(x)<k$, the point $y_{n_k(x)}$ is a {\it unique limit point of $E$}. See \rf{YN}. From this, definition of $a_E$ and \rf{YJ-NX},  for every $j, n_k(x)< j\le k$, we have
$$
a_E(Y_{j-1}(x))=a_E(Y_{n_k(x)}(x))=y_{n_k(x)}=y_j(x),
$$
proving property {\it (a)} in the case under consideration. \bx
\medskip
\begin{lemma}\lbl{NA-P} Let $x_1,x_2\in E$ and let $0\le j\le k$. If $x_1\le x_2$ then
\bel{MN-A}
\min Y_j(x_1)\le \min Y_j(x_2)~~~~~~\text{and}~~~~~~
\max Y_j(x_1)\le \max Y_j(x_2)\,.
\ee
\end{lemma}
\par {\it Proof.} We proceed by induction on $j$. Because $Y_0(x_1)=\{x_1\}$ and $Y_0(x_2)=\{x_2\}$, we conclude that \rf{MN-A} holds for $j=0$.
\par Suppose that \rf{MN-A} holds for some $j, 0\le j\le k-1$. Let us prove that
\bel{MN-J1}
\min Y_{j+1}(x_1)\le \min Y_{j+1}(x_2)\,.
\ee
\par We recall that, thanks to \rf{DF-JY1}, $y_{j+1}(x_\ell)=a_E(Y_j(x_\ell))$ for each $\ell=1,2,$ and
\bel{XL-1}
Y_{j+1}(x_\ell)=Y_{j}(x_\ell)\cup\{y_{j+1}(x_\ell)\}\,.
\ee
\par If $Y_{j}(x_2)$ contains a limit point of $E$, then $y_{j+1}(x_2)\in Y_{j}(x_2)$ so that $Y_{j+1}(x_2)=Y_{j}(x_2)$. This equality and assumption \rf{MN-A} imply that
$$
\min Y_{j+1}(x_1)\le \min Y_{j}(x_1)\le \min Y_{j}(x_2)= \min Y_{j+1}(x_2)
$$
proving \rf{MN-J1} in the case under consideration.
\smallskip
\par Now, suppose that all points of $Y_{j}(x_2)$ are {\it isolated points} of $E$. In particular, from part (b) of Lemma \reff{L-YP} and definitions \rf{YN}, \rf{YJ-NX}, we have $0\le j\le n_k(x_2)$. This inequality and part (c) of Lemma \reff{L-YP} imply that $\#Y_{j}(x_2)=j+1$.
\smsk
\par Consider two cases. First, let us assume that
\bel{BB-J}
\min Y_{j}(x_1)< \min Y_{j}(x_2)\,.
\ee
Then, for each point $a\in\R$ nearest to $Y_{j}(x_2)$ on the set $E\setminus Y_{j}(x_2)$, we have $a\ge \min Y_{j}(x_1)$. This inequality, definition of $a_E$ and \rf{DF-JY1} together yield
$
a_E(Y_{j}(x_2))=y_{j+1}(x_2)\ge \min Y_{j}(x_1)\,.
$
\par Combining this inequality with \rf{XL-1} and \rf{BB-J}, we obtain the required inequality \rf{MN-J1}.
\smallskip
\par Now, prove \rf{MN-J1} whenever $\min Y_{j}(x_1)=\min Y_{j}(x_2)$. This equality and the second inequality in \rf{MN-A} imply that
$$
Y_{j}(x_1)\subset I=[\min Y_{j}(x_2),\max Y_{j}(x_2)]\,.
$$
\par In turn, \rf{MNM-Y} tells us that $I\cap E=Y_{j}(x_2)$ proving that $Y_{j}(x_1)\subset Y_{j}(x_2)$. Recall that in the case under consideration all points of $Y_{j}(x_2)$ are isolated points of $E$. Therefore, all points of $Y_{j}(x_1)$ are isolated points of $E$ as well. Now, using the same argument as for the set $Y_{j}(x_2)$, we conclude that
$\#Y_{j}(x_1)=j+1=\#Y_{j}(x_2).$
\smsk
\par Thus $Y_{j}(x_1)\subset Y_{j}(x_2)$ and $\#Y_{j}(x_1)=\#Y_{j}(x_2)$. Hence, $Y_{j}(x_1)=Y_{j}(x_2)$ proving \rf{MN-J1} in the case under consideration. \smallskip
\par In the same fashion we prove that
$\max Y_{j+1}(x_1)\le \max Y_{j+1}(x_2)$.
\par The proof of the lemma is complete.\bx
\medskip
\begin{lemma}(\cite [p. 231]{KM}) Let $x_1,x_2\in E$, and let $Y_k(x_1)\ne Y_k(x_2)$. Then for all $0\le i,j\le k$ the following inequality
$$
\max\{\,|y_i(x_1)-y_{j}(x_1)|,|y_i(x_2)-y_{j}(x_2)\,|\}\le \max\{i,j\}\,|x_1-x_2|
$$
holds.
\end{lemma}
\par This lemma implies the following
\begin{corollary}\lbl{COR-1} For every $x_1,x_2\in E$ such that $Y_k(x_1)\ne Y_k(x_2)$ the following inequality
$$
\diam Y_k(x_1)+\diam Y_k(x_2)\le 2\,k\,|x_1-x_2|
$$
holds.
\end{corollary}
\bigskip\bigskip
\par {\bf 3.2. Lagrange polynomials and divided differences at interpolation knots.}
\medskip
\addtocontents{toc}{~~~~3.2. Lagrange polynomials and divided differences at interpolation knots. \hfill \thepage\par}
\par In this section we present a series of important properties of Lagrange polynomials which we use later on in proofs of extension criteria.
\begin{lemma}\lbl{DK-T} Let $k$ be a nonnegative integer, and let $P\in\Pc_k$. Suppose that $P$ has $k$ real distinct roots which lie in a set $S\subset\R$. Let $I\subset \R$ be a closed interval.
\par Then for every $i, 0\le i\le k$, the following inequality
$$
\max_I|P^{(i)}|\le (\diam (I\cup S))^{k-i} \,|P^{(k)}|
$$
holds.
\end{lemma}
\par {\it Proof.} Let $x_j$, $j=1,...,k$, be the roots of $P$, and let $X=\{x_1,...,x_k\}$. The lemma's hypothesis tells us that $X\subset S$. Clearly,
$$
P(x)=\frac{P^{(k)}}{k!}\,\prod_{i=1}^k(x-x_i), ~~~x\in\R,
$$
so that for every $i, 0\le i\le k$,
$$
P^{(i)}(x)=\frac{i!}{k!}\,P^{(k)}
\smed_{X'\subset X,\,\,\#X'=k-i}\,\,\pmed_{y\in X'}(x-y),
~~~x\in\R.
$$
Hence,
$$
\max_I |P^{(i)}|\le \frac{i!}{k!}\, \frac{k!}{i!(k-i)!}\,(\diam (I\cup X))^{k-i}\,|P^{(k)}|
\le (\diam (I\cup S))^{k-i}\,|P^{(k)}|
$$
proving the lemma.\bx
\medskip
\par We recall that, given $S\subset\R$ with $\#S=k+1$ and a function $f:S\to\R$, by $L_S[f]$ we denote the Lagrange polynomial of degree at most $k$ interpolating $f$ on $S$.
\begin{lemma}\lbl{LP-T} Let $S_1,S_2\subset \R$, $S_1\ne S_2$, and let $\#S_1=\#S_2=k+1$ where $k$ is a nonnegative integer. Let $I\subset\R$ be a closed interval. Then for every function $f:S_1\cup S_2\to\R$ and every $i, 0\le i\le k$, we have
\bel{LW-11}
\max_I |L_{S_1}^{(i)}[f]-L_{S_2}^{(i)}[f]|\le (k+1)!\,
(\diam (I\cup S_1\cup S_2))^{k-i} \, A
\ee
where
\bel{A-12}
A=\max_{\substack{S'\subset S_1\cup S_2\\\#S'=k+2}}
|\Delta^{k+1}f[S']|\,\diam S'\,.
\ee
\end{lemma}
\par {\it Proof.} Let  $n=k+1-\#(S_1\cap S_2)$; then $n\ge 1$ because $S_1\ne S_2$. Let $\{Y_j:j=0,...,n\}$ be a family of $(k+1)$-point subsets of $S$ such that $Y_0=S_1$, $Y_n=S_2$, and $\#(Y_j\cap Y_{j+1})=k$ for every $j=0,...,n-1$.
\par Let $P_j=L_{Y_j}[f]$, $j=0,...,n$. Then
$$
\max_I |L_{S_1}^{(i)}[f]-L_{S_2}^{(i)}[f]|=
\max_I |P_0^{(i)}-P_n^{(i)}|\le \smed_{j=0}^{n-1}
\max_I |P_j^{(i)}-P_{j+1}^{(i)}|\,.
$$
\par Note that each point $y\in Y_j\cap Y_{j+1}$ is a root of the polynomial $P_j-P_{j+1}\in \Pc_k$. Thus, if the polynomial $P_j-P_{j+1}$ is not identically $0$, it has precisely $k$ distinct real roots which belong to the set $S_1\cup S_2$. We apply Lemma \reff{DK-T}, taking $P=P_j-P_{j+1}$ and $S=S_1\cup S_2$, and obtain the following:
$$
\max_I|P_j^{(i)}-P_{j+1}^{(i)}|\le
(\diam (I\cup S_1\cup S_2))^{k-i} \,|P_j^{(k)}-P_{j+1}^{(k)}|.
$$
From \rf{D-IND} and \rf{D-LAG}, we have
$$
|P_j^{(k)}-P_{j+1}^{(k)}|=
|L^{(k)}_{Y_j}[f]-L^{(k)}_{Y_{j+1}}[f]|=
k!|\Delta^{k}f[Y_j]-\Delta^{k}f[Y_{j+1}]|
\le
k!|\Delta^{k+1}f[S^{(j)}]|\,\diam S^{(j)}
$$
where $S^{(j)}=Y_j\cup Y_{j+1}$, $j=0,...,n-1$. Hence,
\bel{LM-U1}
\max_{I} |L_{S_1}^{(i)}[f]-L_{S_2}^{(i)}[f]|\le k!
\,(\diam (I\cup S_1\cup S_2))^{k-i}\,
\smed_{j=0}^{n-1}\,|\Delta^{k+1}f[S^{(j)}]|\,\diam S^{(j)}.
\ee
\par Clearly, $S^{(j)}\subset S_1\cup S_2$ and $\#S^{(j)}=k+2$. Therefore, each summand of the sum in the right hand side of \rf{LM-U1} is bounded by $A$ (see \rf{A-12}). This, \rf{LM-U1} and inequality $n\le k+1$ imply \rf{LW-11} completing the proof of the lemma.\bx
\begin{lemma}\lbl{LP-TSQ} Let $k$ be a nonnegative integer, $\ell\in\N$, $k<\ell$, and let $\Yc=\{y_j\}_{j=0}^\ell$ be a strictly increasing sequence in $\R$. Let $I=[y_0,y_\ell]$, $S_1=\{y_0,...,y_k\}$, $S_2=\{y_{l-k},...,y_\ell\}$, and let
\bel{SJ-D}
S^{(j)}=\{y_j,...,y_{k+j+1}\},~~~j=0,...,\ell-k-1.
\ee
\par Then for every function $f:\Yc\to\R$, every $i, 0\le i\le k$, and every $p\in[1,\infty)$ the following inequality
\bel{LW-1SQ}
\max_{I} |L_{S_1}^{(i)}[f]-L_{S_2}^{(i)}[f]|^p\le C(k)^p
\,(\diam I)^{(k-i+1)p-1}\,
\smed_{j=0}^{\ell-k-1}\,|\Delta^{k+1}f[S^{(j)}]|^p\,(\diam S^{(j)})
\ee
holds.
\end{lemma}
\par {\it Proof.} Repeating the proof of inequality  \rf{LM-U1}, we obtain the following:
\bel{LM-U}
B=\max_{I} |L_{S_1}^{(i)}[f]-L_{S_2}^{(i)}[f]|\le k!
\,(\diam I)^{k-i}\,
\smed_{j=0}^{\ell-k-1}\,|\Delta^{k+1}f[S^{(j)}]|\,\diam S^{(j)}\,.
\ee
\par Let $I_j=[y_j,y_{k+j+1}]$. Then, thanks to \rf{SJ-D}, $\diam I_j=\diam S_j=y_{k+j+1}-y_j$. Furthermore, since $\{y_j\}_{j=0}^\ell$ is a strictly increasing sequence and $\#S_j=k+2$, the covering multiplicity of the family $\{I_j:j=0,...,\ell-k-1\}$ is bounded by $2k+3$.
Hence,
\bel{L-OU}
\smed_{j=0}^{\ell-k-1}\,\diam S^{(j)}=\smed_{j=0}^{\ell-k-1}\,\diam I_j=\smed_{j=0}^{\ell-k-1}\,|I_j|\le (2k+3)\,|I|\,.
\ee
\par This inequality, H\"{o}lder's inequality and \rf{LM-U} together imply that
$$
B\le
(k!)^p\,(\diam I)^{(k-i)p}\,
\left(\smed_{j=0}^{\ell-k-1}\,\diam S^{(j)}\right)^{p-1}
\,\,\smed_{j=0}^{\ell-k-1}\,|\Delta^{k+1}f[S^{(j)}]|^p
\,\diam S^{(j)}.
$$
\par From this and \rf{L-OU} we have \rf{LW-1SQ} proving the lemma.\bx
\begin{lemma}\lbl{LM-J} Let $k$ be a nonnegative integer and let $1<p<\infty$. Let $f$ be a function defined on a closed set $E\subset\R$ with $\#E>k+1$. Suppose that
\bel{SW-H}
\lambda=\sup_{S\subset E,\,\#S=k+2}\,
|\Delta^{k+1}f[S]|\,(\diam S)^{\frac1p}<\infty\,.
\ee
\par Then for every limit point $x$ of $E$ and every $i,0\le i\le k$, there exists a limit
\bel{FJ}
f_i(x)=\lim_{S\to x,\,S\subset E,\,\#S=k+1}\,\,
L^{(i)}_{S}[f](x)\,.
\ee
(Recall that the notation $S\to x$ means that
$\diam(S\cup \{x\})\to 0$.)
\par Furthermore, let $P_x\in\Pc_k$ be a polynomial such that
\bel{P-XLK}
P^{(i)}_x(x)=f_i(x)~~~\text{for every}~~~i, 0\le i\le k\,.
\ee
Then for every $\delta>0$ and every set $S\subset E$ such that $\#S=k+1$ and $\diam(S\cup \{x\})<\delta$ the following inequality
\bel{M-PL}
\max_{[x-\delta,x+\delta]}
|P^{(i)}_x-L_{S}^{(i)}[f]|\le
C\,\lambda\,\delta^{k+1-i-1/p},~~~0\le i\le k,
\ee
holds. Here $C$ is a constant depending only on $k$.
\end{lemma}
\par {\it Proof.} Let $\delta>0$ and let $S_1,S_2$  be two subsets of $E$ such that $\#S_j=k+1$ and $\diam(S_j\cup\{x\})<\delta$, $j=1,2$. Hence, $S=S_1\cup S_2\subset I=[x-\delta,x+\delta]$. Lemma \reff{LP-T} tells us that
$$
|L_{S_1}^{(i)}[f](x)-L_{S_2}^{(i)}[f](x)|
\le (k+1)!\,(\diam I)^{k-i}\,\max_{S'\subset S,\,\,\#S'=k+2}
|\Delta^{k+1}f[S']|\,\diam S'\,.
$$
Thanks to \rf{SW-H},
$|\Delta^{k+1}f[S']|\le\,\lambda\,(\diam S')^{-\frac1p}
$ for every $(k+2)$-point subset $S'\subset E$, so that
\be
|L_{S_1}^{(i)}[f](x)-L_{S_2}^{(i)}[f](x)|
&\le&
(k+1)!\,\lambda\,(2\delta)^{k-i}\,\max_{S'\subset S,\,\,\#S'=k+2}\,(\diam S')^{1-1/p}\nn\\
&\le&
(k+1)!\,\lambda\,(2\delta)^{k-i}\,(2\delta)^{1-1/p}. \nn
\ee
Hence,
\bel{LP-S}
|L_{S_1}^{(i)}[f](x)-L_{S_2}^{(i)}[f](x)|
\le C(k)\,\lambda\,\delta^{k+1-i-1/p},
\ee
so that
$$
|L_{S_1}^{(i)}[f](x)-L_{S_2}^{(i)}[f](x)|
\to 0~~~\text{as}~~~\delta\to 0\,.
$$
(We recall that $p>1$.) This proves the existence of the limit in \rf{FJ}.
\smsk
\par Let us prove inequality \rf{M-PL}. Thanks to \rf{LP-S}, for every two sets $S, \tS\in E$, with $\#S=\#\tS=k+1$ such that  $\diam(S\cup\{x\}), \diam(\tS\cup\{x\})<\delta$, the following inequality
$$
|L_{\tS}^{(i)}[f](x)-L_{S}^{(i)}[f](x)|
\le C(k)\,\lambda\,\delta^{k+1-i-1/p}
$$
holds. Passing to the limit in this inequality whenever the set $\tS\to x$ (i.e., $\diam(\tS\cup \{x\})\to 0$), we obtain the following:
$$
|P_x^{(i)}(x)-L_{S}^{(i)}[f](x)|
\le C(k)\,\lambda\,\delta^{k+1-i-1/p}\,.
$$
See \rf{FJ} and \rf{P-XLK}. Therefore, for each $y\in[x-\delta,x+\delta]$, we have
\be
|P_x^{(i)}(y)-L_{S}^{(i)}[f](y)|&=&
\left|\smed_{j=0}^{k-i}
\frac{1}{j!}\,
(P_x^{(i+j)}(x)-L_{S}^{(i+j)}[f](x))\,(y-x)^j\right|
\le
\smed_{j=0}^{k-i}
\frac{1}{j!}\,
|P_x^{(i+j)}(x)-L_{S}^{(i+j)}[f](x)|\,\delta^j\nn\\
&\le& C(k)\,\lambda\,\smed_{j=0}^{k-i}\delta^{k+1-i-j-1/p}\,\delta^j
\le C(k)\,\lambda\,\delta^{k+1-i-1/p}\,.\nn
\ee
\par The proof of the lemma is complete.\bx
\begin{lemma}\lbl{D-CN} Let $k,p,E,f,\lambda$ and $x$ be as in Lemma \reff{LM-J}. Then for every $i, 0\le i\le k$,
$$
\lim_{S\to x,\,S\subset E,\,\#S=i+1}\,\,
i!\,\Delta^{i}f[S]=f_i(x)\,.
$$
\end{lemma}
\par {\it Proof.} Let $\delta>0$ and let $S\subset E$ be a finite set such that $\#S=i+1$ and
$\diam(S\cup\{x\})<\delta$. Since $x$ is a limit point of $E$, there exists a set $Y\subset E\cap[x-\delta,x+\delta]$ with $\#Y=k+1$ such that $S\subset Y$. Then, thanks to \rf{M-PL}, 
\bel{M-Y}
\max_{[x-\delta,x+\delta]}
|P^{(i)}_x-L_{Y}^{(i)}[f]|\le
C(k)\,\lambda\,\delta^{k+1-i-1/p}\,.
\ee
\par Because the Lagrange polynomial $L_Y[f]$ interpolates $f$ on $S$, we have $\Delta^{i}f[S]=\Delta^{i}(L_Y[f])[S]$ so that, thanks to \rf{D-KSI}, there exists
$\xi\in [x-\delta,x+\delta]$ such that $i!\,\Delta^{i}f[S]=L_Y^{(i)}[f](\xi)$.
\par This equality and \rf{M-Y} imply that
$$
|P^{(i)}_x(\xi)-i!\Delta^{i}f[S]|=
|P^{(i)}_x(\xi)-L_{Y}^{(i)}[f](\xi)|\le
C(k)\,\lambda\,\delta^{k+1-i-1/p}.
$$
Hence,
\be
|f_i(x)-i!\Delta^{i}f[S]|&=&|P^{(i)}_x(x)-i!\Delta^{i}f[S]|
\le |P^{(i)}_x(x)-P^{(i)}_x(\xi)|+
|P^{(i)}_x(\xi)-i!\Delta^{i}f[S]|\nn\\
&\le&
|P^{(i)}_x(x)-P^{(i)}_x(\xi)|+
C(k)\,\lambda\,\delta^{k+1-i-1/p}.
\nn
\ee
\par Since $P^{(i)}_x$ is a continuous function and $p>1$, the right hand side of this inequality tends to $0$ as $\delta\to 0$ proving the lemma.\bx
\begin{lemma}\lbl{CV-LG} Let $p\in(1,\infty)$, $k\in\N$, and let $f$ be a function defined on a closed set $E\subset\R$ with $\#E>k+1$. Suppose that $f$ satisfies condition \rf{SW-H}.
\par Let $x\in E$ be a limit point of $E$, and let $S$ be a subset of $E$ with $\#S\le k$ containing $x$. Then for every $i, 0\le i\le k+1-\#S$,
$$
\lim_{\substack{S'\setminus S\to x \\
S\subset S'\subset E,\,\#S'=k+1}} \,\,L_{S'}^{(i)}[f](x)=f_i(x)\,.
$$
\end{lemma}
\par {\it Proof.} For $S=\{x\}$ the statement of the lemma follows from Lemma \reff{LM-J}.
\par Suppose that $\#S>1$. Let $I_0=[x-1/2,x+1/2]$ so that $\diam I_0=1$. We prove that for every $i, 0\le i\le k-1$, the family of functions
$$
\{L_{Y}^{(i+1)}[f]: Y\subset I_0\cap E, \#Y=k+1\}
$$
is uniformly bounded on $I_0$ provided condition \rf{SW-H} holds. Indeed, fix a subset $Y_0\subset I_0\cap E$ with $\#Y_0=k+1$. Lemma \reff{LP-T} and \rf{LW-11} together imply that for arbitrary $Y\subset I_0\cap E$, $Y\ne Y_0$, with $\#Y=k+1$, we have
$$
\max_{I_0} |L_{Y}^{(i+1)}[f]-L_{Y_0}^{(i+1)}[f]|\le (k+1)!\,
(\diam (I_0\cup Y\cup Y_0))^{k-i-1} \, A=(k+1)!\,A
$$
where $A=\max\{|\Delta^{k+1}f[S']|\,\diam S':
~S'\subset Y_0\cup Y,~ \#S'=k+2\}$.
\smsk
\par Therefore, thanks to \rf{SW-H},
$$
\max_{I_0} |L_{Y}^{(i+1)}[f]-L_{Y_0}^{(i+1)}[f]|
\le (k+1)! \lambda \max_{S'\subset Y\cup Y_0,\,\#S'=k+2}
(\diam S')^{1-1/p}\le (k+1)! \lambda\,.
$$
\par Applying this inequality to an arbitrary set $Y\subset I_0\cap E$ with $\#Y=k+1$ and to every $i, 0\le i\le k-1$, we conclude that
\bel{UB-L}
\max_{I_0} |L_{Y}^{(i+1)}[f]|\le B_i~~~~~~\text{where}~~~~~~
B_i=\max_{I_0}|L_{Y_0}^{(i+1)}[f]|+(k+1)! \lambda.
\ee
\par Fix $\ve>0$. Lemma \reff{D-CN} tells us that there exists $\dw\in(0,1/2]$ such that for an arbitrary set $V\subset E$, with $\diam (V\cup\{x\})<\dw$ and $\#V=i+1$, the following inequality
\bel{EP-1}
|i!\,\Delta^{i}f[V]-f_i(x)|\le \ve/2
\ee
holds. Let $S'$ be an arbitrary subset of $E$ such that $S\subset S'$, $\#S'=k+1$, and
\bel{BI-D}
\diam ((S'\setminus S)\cup\{x\}) <\delta=\min\{\dw,\ve/(2B_i)\}\,.
\ee
\par Recall that $\#S'-\#S=k+1-\#S\ge i$, so that there exists a subset $V\subset (S'\setminus S)\cup\{x\}$ with $\#V=i+1$. Thanks to \rf{D-KSI}, there exists $\xi\in[x-\delta,x+\delta]$ such that
$i!\,\Delta^i(L_{S'}[f])[V]=L^{(i)}_{S'}[f](\xi)$.
\smsk
\par On the other hand, because the polynomial $L_{S'}[f]$ interpolates $f$ on $V$, the divided difference $\Delta^if[V]=\Delta^i(L_{S'}[f])[V]$ so that, $i!\Delta^if[V]=L^{(i)}_{S'}[f](\xi)$.
\par This and \rf{EP-1} imply that $|f_i(x)-L^{(i)}_{S'}[f](\xi)|\le \ve/2$.
\smsk
\par It remains to note that, thanks to \rf{UB-L} and \rf{BI-D},
$$
|L^{(i)}_{S'}[f](\xi)-L^{(i)}_{S'}[f](x)|\le
\left(\max_{[x-1/2,x+1/2]}|L^{(i+1)}_{S'}[f]|\right)\cdot |x-\xi|\le B_i\, \delta\le B_i\,(\ve/(2B_i))=\ve/2,
$$
so that
$$
|f_i(x)-L^{(i)}_{S'}[f](x)|\le
|f_i(x)-L^{(i)}_{S'}[f](\xi)|+
|L^{(i)}_{S'}[f](\xi)-L^{(i)}_{S'}[f](x)|\le
\ve/2+\ve/2=\ve
$$
proving the lemma.\bx
\medskip\bigskip
\par {\bf 3.3. Whitney $m$-fields and Hermite polynomials.}
\medskip
\addtocontents{toc}{~~~~3.3. Whitney $m$-fields and Hermite polynomials. \hfill \thepage\par}

\par We turn to constructing of the Whitney $(m-1)$-field
$\VP^{(m,E)}$ mentioned at the beginning of Section 3. See \rf{PME-1}. Everywhere in this section we will assume that $f$ is a function on $E$ satisfying the following condition:
\bel{A-FE}
\sup_{S\subset E,\,\#S=m+1}\,
|\Delta^{m}f[S]|\,(\diam S)^{\frac1p}<\infty\,.
\ee
\par Let $k=m-1$. Given $x\in E$, let
\bel{S-X-D}
\SH_x=Y_{k}(x)=\{y_0(x),...,y_{n_k(x)}(x)\}
\ee
and let
\bel{S-XSM}
s_x=y_{n_k(x)}\,.
\ee
\par We recall that the points $y_j(x)$ and the sets $Y_j(x)$ are defined by formulae \rf{DF-JY1}-\rf{YJ-NX}.
\medskip
\par The next two propositions describe the main properties of the sets $\{\SH_x: x\in E\}$ and the points $\{s_x: x\in E\}$. These properties are immediate from Lemmas \reff{L-YP}, \reff{NA-P} and Corollary \reff{COR-1}.
\begin{proposition}\lbl{SET-SX}
\par (i) $\SH_x\subset E$,\, $x\in \SH_x$ and~ $\#\SH_x\le m$~ for every $x\in E$. Furthermore,
\bel{MNM-X1}
[\min \SH_x\,,\max \SH_x]\cap E= \SH_x\,;
\ee
\par (ii) For every $x_1,x_2\in E$ such that $\SH_{x_1}\ne \SH_{x_2}$ the following inequality
\bel{D-F12}
\diam \SH_{x_1}+\diam \SH_{x_2}\le 2\,m\,|x_1-x_2|
\ee
holds;
\smallskip 
\par (iii). If $x_1,x_2\in E$ and $x_1<x_2$ then
$$
\min \SH_{x_1}\le \min \SH_{x_2}~~~\text{and}~~~
\max \SH_{x_1}\le \max \SH_{x_2}.
$$
\end{proposition}
\begin{proposition}\lbl{PR-SX}
\par (i) The point
\bel{SX-S}
s_x~~~\text{belongs to}~~~\SH_x
\ee
for every $x\in E$. This point is either minimal or maximal point of the set \,$\SH_x$;
\smallskip
\par (ii) All points of the set $\SH_x\setminus \{s_x\}$ are isolated points of $E$ provided $\#S_x>1$. If $y\in \SH_x$ and $y$ is a limit point of $E$, then $y=s_x$;
\smsk
\par (iii) If~ $\#\SH_x<m$ then $s_x$ is a limit point of $E$.
\end{proposition}
\begin{definition}\lbl{P-X} {\em \par Given a function $f:E\to\R$ satisfying condition \rf{A-FE}, we define the Whitney $(m-1)$-field~ $\VP^{(m,E)}[f]=\{P_x\in\PMO:x\in E\}$\, as follows:
\medskip
\par (i) If~ $\#\SH_x<m$, part (iii) of Proposition \reff{PR-SX} tells us that $s_x$ is a {\it limit point} of $E$. Then, thanks to \rf{A-FE} and  Lemma \reff{LM-J}, for every $i, 0\le i\le m-1$, there exists a limit
\bel{FJ-H}
f_i(s_x)=\lim_{\substack{S\to s_x\\
S\subset E,\,\#S=m}}\,\,
L^{(i)}_{S}[f](s_x)\,.
\ee
\par We define a polynomial $P_x\in\PMO$ as the Hermite polynomial satisfying the following conditions:
\bel{PX-IS}
P_x(y)=f(y)~~~\text{for every}~~~y\in \SH_x,               \ee
and
\bel{PX-SX}
P^{(i)}_x(s_x)=f_i(s_x)~~~\text{for every}~~~i, 1\le i\le m-\#\SH_x\,.
\ee
\par (ii) If~ $\#\SH_x=m$, we put
\bel{PX-M}
P_x=L_{\SH_x}[f].
\ee
}
\end{definition}
\par The next lemma shows that the field $\VP^{(m,E)}[f]$ determined by Definition \reff{P-X} is well defined.
\begin{lemma}\lbl{P-WDEF} For each $x\in E$ there exists the unique polynomial $P_x$ satisfying conditions \rf{PX-IS} and \rf{PX-SX} provided condition \rf{A-FE} holds.
\end{lemma}
\par {\it Proof.} In case (i) ($\#\SH_x<m$) the existence and uniqueness of $P_x$ satisfying \rf{PX-IS} and \rf{PX-SX} is immediate from \cite[Ch. 2, Section 11]{BZ}. See also formula \rf{P-HPR} below.
\par Clearly, in case (ii) ($\#\SH_x=m$) property \rf{PX-IS} holds as well, and \rf{PX-SX} holds vacuously.
\bx
\smsk
\par Let us note that $x\in \SH_x$ and
$P_x=f$ on $\SH_x$ (see \rf{PX-IS}) proving that
\bel{PX-X1}
P_x(x)=f(x)~~~\text{for every}~~~x\in E\,.
\ee
\par For the case $m>1$ and $\#\SH_x<m$, we give an explicit formula for the Hermite polynomials $P_x, x\in E$, from Definition \reff{P-X}. See \cite[Ch. 2, Section 11]{BZ}.
\par Let $n=\#\SH_x-1$ and let $y_i=y_i(x)$, $i=0,...,n$, so that $\SH_x=\{y_0,...,y_n\}$. See \rf{S-X-D}. (Note also that in these settings $s_x=y_n$.) In this case the Hermite polynomial $P_x$ satisfying \rf{PX-IS} and \rf{PX-SX} can be represented as a linear combination of polynomials
$$
H_0,..., H_n, \tH_1,....,\tH_{m-n-1}\in \PMO
$$
which are uniquely determined by the following conditions:
\bigskip
\par (i) $H_i(y_i)=1$ for every $i,~ 0\le i\le n,$~~ and ~~$H_j(y_i)=0$~ for every~ $0\le i,j\le n,~i\ne j$, and
$$
H'_i(y_n)=...=H_i^{(m-n-1)}(y_n)=0~~~\text{for every}~~
i,~0 \le i\le n\,.
$$
\par (ii) $\tH_j(y_i)=0$~ for every $0\le i\le n, 1\le j\le m-n-1,$ and for every $1\le j\le m-n-1$,
$$
\tH_j^{(j)}(y_n)=1 ~~~~\text{and}~~~~ \tH^{(\ell)}_j(y_n)=0~~~\text{for every}~~ \ell,~1\le \ell\le m-n-1,~\ell\ne j\,.
$$
\par The existence and uniqueness of $H_i$ and $\tH_j$, $0\le i\le n$, $1\le j\le m-n-1$, are proven in  \cite[Ch. 2, Section 11]{BZ}. It is also shown there that for every $P\in\PMO$ the following unique representation
\bel{P-BS}
P(y)=\smed_{i=0}^n\,P(y_i)\,H_i(y)\,+\,
\smed_{j=1}^{m-n-1}\,P^{(j)}(y_n)\,\tH_j(y),~~~~~~y\in\R,
\ee
holds. In particular,
\bel{P-HPR}
P_x(y)=\smed_{i=0}^n\,f(y_i)\,H_i(y)\,+\,
\smed_{j=1}^{m-n-1}\,f_j(y_n)\,\tH_j(y),~~~y\in\R\,.
\ee
Clearly, $P_x$ meets conditions \rf{PX-IS} and \rf{PX-SX}.
\medskip
\par Let $I\subset \R$ be a bounded closed interval, and let $C^m(I)$ be the space of all $m$-times continuously differentiable functions on $I$. We norm $C^m(I)$ by
$$
\|f\|_{C^m(I)}=\smed_{i=0}^m\max_I|f^{(i)}|\,.
$$
\par We will need the following important property of the polynomials $\{P_x: x\in E\}$.
\begin{lemma}\lbl{HP-CN} Let $f$ be a function defined on a closed set $E\subset\R$ with $\#E>m+1$, and satisfying condition \rf{A-FE}. Let $I$ be a bounded closed interval in $\R$. Then for every $x\in E$
$$
\lim_{\substack{S'\setminus \SH_x\to s_x \smallskip\\
\SH_x\subset\, S'\subset E,\,\,\#S'=m}} \,\,\|L_{S'}[f]-P_x\|_{C^m(I)}=0\,.
$$
\end{lemma}
\par {\it Proof.} The lemma is obvious whenever $\# \SH_x=m$ because in this case  $L_{\SH_x}[f]=P_x$. In particular, the lemma is trivial for $m=1$.
\par Let now $m>1$ and let $\# \SH_x<m$. In this case $P_x$ can be represented in the form \rf{P-HPR}. Because $s_x$ is a limit point of $E$ (see part (iii) of Proposition \reff{PR-SX}), Lemma \reff{CV-LG} and \rf{PX-SX} imply that
\bel{L-YNC}
\lim_{\substack{S'\setminus \SH_x\to s_x \smallskip\\
\SH_x\subset\, S'\subset E,\,\,\#S'=m}} L_{S'}^{(i)}[f](s_x)=f_i(s_x)=P_x^{(i)}(s_x)
~~~\text{for every}~~~1\le i\le m-n-1\,.
\ee
\par Let $n=\#\SH_x-1$ and let
$\SH_x=\{y_0,...,y_n\}$ where $y_i=y_i(x)$,~ $i=0,...,n$.
\par Then, thanks to \rf{P-BS}, for every set $S'\subset E$ with $\#S'=m $ such that $\SH_x\subset\, S'$, the polynomial $L_{S'}[f]$ has the following representation:
$$
L_{S'}[f](y)=\smed_{i=0}^n\,f(y_i)\,H_i(y)\,+\,
\smed_{j=1}^{m-n-1}\,L^{(j)}_{S'}[f](s_x)\,\tH_j(y),
~~~y\in\R\,.
$$
\par From this, \rf{P-HPR} and \rf{L-YNC}, we have
$$
\max_I\,\left|L_{S'}[f]-P_x\right|=\max_I
\,\left|\smed_{j=1}^{m-n-1}\,
(L^{(j)}_{S'}[f](s_x)-P^{(j)}_x(s_x))\,\tH_j\right|
\le
\,\smed_{j=1}^{m-n-1}\,
\left|L^{(j)}_{S'}[f](s_x)-P^{(j)}_x(s_x)\right|
\,\max_I\left|\tH_j\right|.
$$
Hence, $\max_I\,\left|L_{S'}[f]-P_x\right|\to 0$
as $S'\setminus \SH_x\to s_x$ provided $\SH_x\subset\, S'\subset E$ and $\#S'=m$.
\par Note that the uniform norm on $I$ and the $C^m(I)$-norm are equivalent norms on the finite dimensional space $\PM$. Therefore, convergence of $L_{S'}[f]$ to $P_x$ in the uniform norm on $I$ implies convergence of $L_{S'}[f]$ to $P_x$ in the $C^m(I)$-norm, proving the lemma.\bx

\bigskip
\par {\bf 3.4. Extension criteria in terms of sharp maximal functions: sufficiency.}
\medskip
\addtocontents{toc}{~~~~3.4. Extension criteria in terms of sharp maximal functions: sufficiency. \hfill \thepage\VST\par}
\par Let $f$ be a function on $E$ such that
$\SHF\in \LPR$. See \rf{SH-F} and \rf{SH-F-EQ}.
\smallskip
\par Let us prove that $f$ satisfies condition \rf{A-FE}. Indeed, let $S=\{x_0,...,x_m\}\subset E$, $x_0<...<x_m$. Clearly, for every $x\in [x_0,x_m]$, we have
$\diam (\{x\}\cup S)=\diam S=x_m-x_0$, so that, thanks to \rf{SH-F-EQ},
$$
|\Delta^mf[S]|^p\,\diam S=
\frac{|\Delta^mf[S]|^p (\diam S)^p}
{(\diam (\{x\}\cup S))^p}\,(x_m-x_0)
\le 2^p\,(\SHF(x))^p\,(x_m-x_0).
$$
\par Integrating this inequality (with respect to $x$) over the interval $[x_0,x_m]$, we obtain the following:
$$
|\Delta^mf[S]|^p\,\diam S
\le 2^p\,\intl_{x_0}^{x_m}\,(\SHF(x))^p\,dx\le 2^p\, \|\SHF\|^p_{\LPR}.
$$
Hence,
$$
\sup_{S\subset E,\,\#S=m+1}\,
|\Delta^{m}f[S]|\,(\diam S)^{\frac1p}
\le 2\,\|\SHF\|_{\LPR}<\infty
$$
proving \rf{A-FE}.
\smallskip
\par This condition and Lemma \reff{P-WDEF} guarantee that the Whitney $(m-1)$-field $\VP^{(m,E)}[f]$ from Definition \reff{P-X} is well defined.
\smsk
\par Now, let us to show that inequality \rf{TR-PW} holds. Its proof relies on Theorem \reff{JET-S} below which provides a criterion for the restrictions of Sobolev jets.
\par For each family $\VP=\{P_x\in\PMO: x\in E\}$ of polynomials we let $\VSH$ denote a certain kind of a ``sharp maximal function'' associated with $\VP$ which is defined by
$$
\VSH(x)=\sup_{a_1,\,a_2\in E,\,\, a_1\ne a_2}\,\, \frac{|P_{a_1}(x)-P_{a_2}(x)|}
{|x-a_1|^{m}+|x-a_2|^{m}},~~~~~~~x\in\R.
$$
\begin{theorem} \lbl{JET-S}(\cite{Sh5}) Let $m\in\N$, $p\in(1,\infty)$, and let $E$ be a closed subset of $\R$. Suppose we are given a family $\VP=\{P_x: x\in E\}$ of polynomials of degree at most $m-1$ indexed by points of $E$.
\par Then there exists a $C^{m-1}$-function $F\in\LMPR$ such that $T_{x}^{m-1}[F]=P_{x}$ \ for every $x\in E$
if and only if\, $\VSH\in\LPR$. Furthermore,
\bel{PME-VS}
\PME\sim \|\VSH\|_{\LPR}
\ee
with the constants in this equivalence depending only on $m$ and $p$.
\end{theorem}
\par We recall that the quantity $\PME$ is defined by \rf{N-VP}.
\begin{lemma}\lbl{PE-1} Let $f$ be a function on $E$ such that $\SHF\in \LPR$. Then for every $x\in\R$ the following inequality
\bel{SMF-C}
(\VP^{(m,E)}[f])^\sharp_{m,E}(x)
\le C(m)\,\SHF(x)~~~
\ee
holds.
\end{lemma}
\par {\it Proof.} Let $x\in\R$, $a_1,a_2\in E$, $a_1\ne a_2$, and let $\rl= |x-a_1|+|x-a_2|$. Let $\TSH_j=S_{a_j}$ and let $s_j=s_{a_j}$, $j=1,2$. See \rf{S-X-D} and \rf{S-XSM}. We know that $a_j,s_{j}\in \TSH_j$, $j=1,2$ (see Propositions \reff{SET-SX} and \reff{PR-SX}).
\par Suppose that $\TSH_1\ne \TSH_2$. From inequality \rf{D-F12}, we have
\bel{DM-S}
\diam \TSH_1+\diam \TSH_2\le 2\,m\,|a_1-a_2|\,.
\ee
\par Fix an $\ve>0$. Lemma \reff{HP-CN} tells us that for each $j=1,2$ there exists an $m$-point subset $\SH_j\subset E$, such that ${\tS\hspace*{-0.7mm}}_j\subset \SH_j$,~
$\diam (\{s_j\}\cup (\SH_j\setminus \TSH_j))\le r$ and
and
\bel{J-2}
|P_{a_j}(x)-L_{\SH_j}[f](x)|\le \ve\,\rl^{m}/2^{m+1}\,.
\ee
\par Recall that $s_j\in\TSH_j$, $j=1,2$, so that
$$
\diam \SH_j\le \diam \TSH_j+\diam(\{s_j\}\cup (\SH_j\setminus \TSH_j))\le\diam \TSH_j+r.
$$
This inequality together with \rf{DM-S} imply that
\bel{J-1}
\diam \SH_j\le 2m\,|a_1-a_2|+r, ~~~j=1,2.
\ee
\par Let $I$ be the smallest closed interval containing $S_1\cup S_2\cup\{x\}$. Since $a_j\in\TSH_j\subset \SH_j$, $j=1,2$,
$$
\diam I\le |x-a_1|+|x-a_2|+\diam \SH_1+\diam \SH_2
$$
so that, thanks to \rf{J-1},
\bel{DI-E}
\diam I\le |x-a_1|+|x-a_2|+4m|a_1-a_2|+2r\le (4m+3)\,r.
\ee
(Recall that $r= |x-a_1|+|x-a_2|$.) From this inequality and inequality \rf{J-2}, we have
\be
|P_{a_1}(x)-P_{a_2}(x)|&\le& |P_{a_1}(x)-L_{\SH_1}[f](x)|+
|L_{\SH_1}[f](x)-L_{\SH_2}[f](x)|
\nn\\
&+&
|P_{a_2}(x)-L_{\SH_2}[f](x)|=J+\ve\rl^{m}/2^{m}
\nn
\ee
where $J=|L_{\SH_1}[f](x)-L_{\SH_2}[f](x)|$.
\smallskip
\par Let us estimate $J$. We may assume that $\SH_1\ne \SH_2$; otherwise $J=0$. We apply Lemma \reff{LP-T} taking $k=m-1$ and $i=0$, and get
$$
J\le\max_{I} |L_{\SH_1}[f]-L_{\SH_2}[f]|\le
m!\,(\diam I)^{m-1}\max_{S'\subset S,\,\,\#S'=m+1}
|\Delta^{m}f[S']|\,\diam S'
$$
where $S=\SH_1\cup \SH_2$. This inequality together with  \rf{DI-E} and \rf{SH-F-EQ} implies that
$$
J\le
C(m)\,r^{m-1}\,(\Delta^m f)^\sharp_E(x)
\max_{S'\subset S,\,\,\#S'=m+1}
\,\diam(\{x\}\cup S')\le C(m)\,r^{m}\,(\Delta^m f)^\sharp_E(x)\,.
$$
\par We are in a position to prove inequality \rf{SMF-C}. We have:
\be
\frac{|P_{a_1}(x)-P_{a_2}(x)|}
{|x-a_1|^{m}+|x-a_2|^{m}}&\le&
2^m\,\rl^{-m}\,|P_{a_1}(x)-P_{a_2}(x)|\le 2^m\,\rl^{-m}\,(J+\ve\rl^{m}/2^{m})\nn\\
&\le&
C(m)2^m\,\rl^{-m}\,r^{m}\,(\Delta^m f)^\sharp_E(x)+\ve
=C(m)\,(\Delta^m f)^\sharp_E(x)+\ve\,.
\nn
\ee
Because $\ve>0$ is arbitrary, we conclude that
\bel{R-Y}
\frac{|P_{a_1}(x)-P_{a_2}(x)|}
{|x-a_1|^{m}+|x-a_2|^{m}}\le
\,C(m)\,(\Delta^m f)^\sharp_E(x)
\ee
provided $\TSH_1\ne \TSH_2$. Clearly, this inequality also holds whenever $\TSH_1=\TSH_2$ because in this case $P_{a_1}=P_{a_2}$.
\par Finally, taking the supremum in the left hand side of \rf{R-Y} over all $a_1,\,a_2\in E$, $a_1\ne a_2$, we obtain \rf{SMF-C}. The proof of the lemma is complete.\bx
\medskip
\par We finish the proof of Theorem \reff{R1-CR} as follows.
\par Let $f$ be a function on $E$ such that $\SHF\in \LPR$, and let $\VP^{(m,E)}[f]=\{P_x\in\PMO:x\in E\}$ be the Whitney $(m-1)$-field from Definition \reff{P-X}. Lemma \reff{PE-1} tells us that
$$
\|(\VP^{(m,E)}[f])^\sharp_{m,E}\|_{\LPR}
\le C(m)\,\|\SHF\|_{\LPR}\,.
$$
Combining this inequality with equivalence \rf{PME-VS}, we obtain inequality \rf{TR-PW}.
\par This inequality and definition \rf{N-VP} imply the existence of a function $F\in\LMPR$ such that $T^{m-1}_x[F]=P_x$ on $E$ and
\bel{T-F}
\|F\|_{\LMPR}\le
2\,\|\VP^{(m,E)}[f]\|_{m,p,E}\le C(m,p)\,\|\SHF\|_{\LPR}\,.
\ee
We also note that $P_x(x)=f(x)$ on $E$, see \rf{PX-X1}, so that
$$
F(x)=T^{m-1}_x[F](x)=P_x(x)=f(x),~~~x\in E\,.
$$
\par Thus $F\in\LMPR$ and $F|_E=f$ proving that $f\in \LMPR|_E$. Furthermore, thanks to \rf{N-LMPR} and \rf{T-F},
$$
\|f\|_{\LMPR|_E}\le \|F\|_{\LMPR}\le C(m,p)\,\|\SHF\|_{\LPR}.
$$
\par The proof of Theorem \reff{R1-CR} is complete.\bx

\SECT{4. A variational extension criterion for Sobolev jets.}{4}
\addtocontents{toc}{4. A variational extension criterion for Sobolev jets. \hfill \thepage\par\VST}

\indent\par One of the main ingredients of our proof of the sufficiency part of Theorem \reff{MAIN-TH} (see Section 6) is the following refinement of Theorem \reff{JET-S}.
\begin{theorem} \lbl{JET-V} Let $m\in\N$, $p\in(1,\infty)$, and let $E\subset\R$ be a closed set. Suppose we are given a Whitney $(m-1)$-field $\VP=\{P_x: x\in E\}$ defined on $E$. There exists a $C^{m-1}$-function $F\in\LMPR$ such that
\bel{TX-P}
T_{x}^{m-1}[F]=P_{x}~~~\text{for every}~~~x\in E
\ee
if and only if the following quantity
\bel{VP-V}
\Nc_{m,p,E}(\VP)=\sup\left\{\,\smed_{j=1}^{k-1}\,\,
\smed_{i=0}^{m-1}\,\,
\frac{|P^{(i)}_{x_j}(x_j)-P^{(i)}_{x_{j+1}} (x_j)|^p}
{(x_{j+1}-x_{j})^{(m-i)p-1}}\right\}^{1/p}
\ee
is finite. Here the supremum is taken over all integers $k>1$ and all finite strictly increasing sequences $\{x_1,...,x_k\}\subset E$.
\par Furthermore, $\PME\sim \Nc_{m,p,E}(\VP)$ with the
constants in this equivalence depending only on $m$.
\end{theorem}
\par {\it Proof.} {\it (Necessity.)} Let $\{x_j\}_{j=1}^k$ be a strictly increasing sequence in $E$. Let $\VP=\{P_x: x\in E\}$ be a Whitney $(m-1)$-field on $E$, and let $F\in\LMPR$ be a function satisfying condition \rf{TX-P}. The Taylor formula with the reminder in the integral form tells us that for every $x\in\R$ and every $a\in E$ the following equality
$$
F(x)-T_{a}^{m-1}[F](x)=\frac{1}{(m-1)!}\intl_a^x\, F^{(m)}(t)\,(x-t)^{m-1}\,dt
$$
holds.
\par Let $0\le i\le m-1$. Differentiating this equality $i$ times (with respect to $x$) we obtain the following:
$$
F^{(i)}(x)-(T_{a}^{m-1}[F])^{(i)}(x)
=\frac{1}{(m-1-i)!}\intl_a^x\, F^{(m)}(t)\,(x-t)^{m-1-i}\,dt.
$$
From this and \rf{TX-P}, we have
$$
P_x^{(i)}(x)-P_{a}^{(i)}(x)=\frac{1}{(m-1-i)!}\intl_a^x\, F^{(m)}(t)\,(x-t)^{m-1-i}\,dt~~~~\text{for every}~~~x\in E.
$$
\par Therefore,  for every $j\in\{1,...,k-1\}$ the following equality
$$
P_{x_j}^{(i)}(x_j)-P_{x_{j+1}}^{(i)}(x_{j})=
\frac{1}{(m-1-i)!}\intl_{x_{j+1}}^{x_{j}}\, F^{(m)}(t)\,(x_{j}-t)^{m-1-i}\,dt
$$
holds. Hence,
$$
\left|P_{x_j}^{(i)}(x_j)-P_{x_{j+1}}^{(i)}(x_{j})\right|^p
\le
\frac{(x_{j+1}-x_j)^{(m-1-i)p}}{((m-1-i)!)^p}
\,\left(\intl_{x_j}^{x_{j+1}}\,|F^{(m)}(t)|\,dt\right)^p
$$
proving that
$$
\frac{|P^{(i)}_{x_j}(x_j)-P^{(i)}_{x_{j+1}} (x_j)|^p}
{(x_{j+1}-x_{j})^{(m-i)p-1}}\le
\frac{(x_{j+1}-x_j)^{1-p}}{((m-1-i)!)^p}
\,\left(\intl_{x_j}^{x_{j+1}}\,|F^{(m)}(t)|\,dt\right)^p
\le \frac{1}{((m-1-i)!)^p}
\intl_{x_j}^{x_{j+1}}\,|F^{(m)}(t)|^p\,dt.
$$
\par Consequently,
\be
\smed_{j=1}^{k-1}\,\,\smed_{i=0}^{m-1}\,\,
\frac{|P^{(i)}_{x_j}(x_j)-P^{(i)}_{x_{j+1}} (x_j)|^p}
{(x_{j+1}-x_{j})^{(m-i)p-1}}&\le&
\smed_{j=1}^{k-1}\,\,\smed_{i=0}^{m-1}
\frac{1}{((m-1-i)!)^p}
\intl_{x_j}^{x_{j+1}}\,|F^{(m)}(t)|^p\,dt
\nn\\
&=&
\left(\smed_{i=0}^{m-1}
\frac{1}{((m-1-i)!)^p}\right)\cdot
\intl_{x_1}^{x_{k}}\,|F^{(m)}(t)|^p\,dt
\le \,e^p\,\|F\|_{\LMPR}^p.
\nn
\ee
\par Taking the supremum in the left hand side of this inequality over all finite strictly increasing sequences $\{x_j\}_{j=1}^k$ in $E$, and then the infimum in the right hand side over all function $F\in\LMPR$ satisfying \rf{TX-P}, we obtain the required inequality
$\Nc_{m,p,E}(\VP)\le e\,\PME$.
\par The proof of the necessity is complete.\bx
\medskip
\par {\it (Sufficiency.)} Let $\VP=\{P_x: x\in E\}$ be
a Whitney $(m-1)$-field defined on $E$ such that
$$
\lambda=\Nc_{m,p,E}(\VP)<\infty.
$$
See \rf{VP-V}. Thus, for every strictly increasing sequence $\{x_j\}_{j=1}^k$ in $E$ the following inequality
\bel{L-N}
\smed_{j=1}^{k-1}\,\,
\smed_{i=0}^{m-1}\,\,
\frac{|P^{(i)}_{x_j}(x_j)-P^{(i)}_{x_{j+1}} (x_j)|^p}
{(x_{j+1}-x_{j})^{(m-i)p-1}}\le \lambda^p
\ee
holds. Our aim is to prove the existence of a function $F\in\LMPR$ such that
$T_{x}^{m-1}[F]=P_{x}$ for every $x\in E$ and $\|F\|_{\LMPR}\le C(m)\,\lambda$.
\msk
\par We construct $F$ with the help of the classical Whitney extension method \cite{W1}. It is proven in \cite{Sh5} that this method provides an almost optimal extension of the restrictions of Whitney $(m-1)$-fields generated by Sobolev $\WMP$-functions. In this paper we will use a special one dimensional version of this method suggested by Whitney in \cite[Section 4]{W2} .
\par Because $E$ is a closed subset of $\R$, the complement of $E$, the set $\R\setminus E$, can be represented as a union of a certain finite or countable family
\bel{TA-E}
\Jc_E=\{J_k=(a_k,b_k): k\in \Kc\}
\ee
of pairwise disjoint open intervals (bounded or unbounded). Thus, $a_k,b_k\in E\cup\{\pm\infty\}$ for all $k\in\Kc$,
$$
\R\setminus E= \mcup\{J_k=(a_k,b_k): k\in\Kc\}~~~\text{and}~~~J_{k'}\mcap\, J_{k''}=\emp~~
\text{for every}~~k',k''\in \Kc, k'\ne k''.
$$
\par To each interval $J\in\Jc_E$ we assign a polynomial $H_J\in \Pc_{2m-1}$ as follows:\smallskip
\par (1)~ Let $J=(a,b)$ be an {\it unbounded} open interval, i.e., either $a=-\infty$ and $b$ is finite, or $a$ is finite and $b=+\infty$. In the first case (i.e., $J=(a,b)=(-\infty,b)$) we set $H_J=P_b$, while in the second case (i.e., $J=(a,b)=(a,+\infty)$) we set
$H_J=P_a$.
\smallskip
\par (2)~ Let $J=(a,b)\in\Jc_E$ be a {\it bounded} interval so that  $a,b\in E$. In this case we define the polynomial $H_J\in \Pc_{2m-1}$ as the Hermite polynomial satisfying the following conditions:
\bel{H-J}
H^{(i)}_J(a)=P^{(i)}_a(a)~~~\text{and}~~~
H^{(i)}_J(b)=P^{(i)}_b(b)~~~\text{for all}~~~
i=0,...,m-1.
\ee
\par The existence and uniqueness of the polynomial $H_J$ follows from \cite[Ch. 2, Section 11]{BZ}.
\smsk
\par Finally, we define the extension $F$ by the formula:
\bel{DEF-F}
F(x)=\left \{
\begin{array}{ll}
P_x(x),& x\in E,\smallskip\\
\sbig\limits_{J\in \Jc_E}\,\,
H_J(x)\,\chi_J(x),& x\in\R\setminus E.
\end{array}
\right.
\ee
\smsk
\par We note that inequality \rf{L-N} implies the following property of the Whitney field $\VP=\{P_x: x\in E\}$: for every $x,y\in E$ and every $i, 0\le i\le m-1$, we have
$$
|P^{(i)}_x(x)-P^{(i)}_y(x)|
\le \lambda\,|x-y|^{m-i-1/p}.
$$
Recall that $p>1$, which implies that
\bel{O-W}
P^{(i)}_x(x)-P^{(i)}_y(x)=o(|x-y|^{m-1-i})
~~~\text{provided}~~~x,y\in E~~~\text{and}~~~ 0\le i\le m-1.
\ee
\par Whitney \cite{W2} proved that for every $(m-1)$-field $\VP=\{P_x: x\in E\}$ satisfying \rf{O-W}, the extension $F$ defined by formula \rf{DEF-F} is a $C^{m-1}$-function on $\R$ which agrees with $\VP$ on $E$, i.e.,
\bel{F-AG-P}
F^{(i)}(x)=P^{(i)}_x(x)~~~~\text{for all}~~~~x\in E
~~~\text{and}~~~i=0,...m-1.
\ee
\par Let us show that $F\in\LMPR$ and $\|F\|_{\LMPR}\le C(m)\,\lambda$. Our proof of these properties of the function $F$ relies on the following description of $\LOPR$-functions.
\begin{theorem}\lbl{CR-SOB} Let $p>1$ and let $\tau>0$. Let $G$ be a continuous function on $\R$ satisfying the following condition: There exists a constant $A>0$ such that for every finite family $\Ic=\{I=[u_I,v_I]\}$ of pairwise disjoint closed intervals of diameter at most
$\tau$ the following inequality
$$
\smed_{I=[u_I,v_I]\in\Ic}\,\frac{|\,G(u_I)-G(v_I)|^p}
{(v_I-u_I)^{p-1}} \le A
$$
holds. Then $G\in\LOPR$ and
$\|G\|_{\LOPR}\le C\,A^{\frac1p}$ where $C$ is an absolute constant.
\end{theorem}
\par {\it Proof.} The Riesz theorem \cite{R} tells us that $G\in\LOPR$. See also \cite{JS}. For $\tau=\infty$ inequality $\|G\|_{\LOPR}\le C\,A^{\frac1p}$ follows from
\cite[Theorem 2]{B} and \cite[Theorem 4]{B2}. (See also
a description of Sobolev spaces obtained in \cite[\S\,4, $3^{\circ}$]{B2}.) For the case $0<\tau<\infty$ we refer the reader to \cite[Section 7]{Sh5}.\bx
\smallskip
\par We will also need the following auxiliary lemmas.
\begin{lemma}\lbl{DF-N} Let $J=(a,b)\in\Jc_E$ be a bounded interval, and let $H_J\in\Pc_{2m-1}$ be the Hermite polynomial satisfying \rf{H-J}. Then for every $n\in\{0,...,m\}$ and every $x\in [a,b]$ the following inequality
$$
|H_J^{(n)}(x)|\le C(m)\,\min\left\{Y_1(x),Y_2(x)\right\}
$$
holds. Here
$$
Y_1(x)=|P_a^{(n)}(x)|+
\left\{\smed_{i=0}^{m-1}\,\,
\frac{|P^{(i)}_b(b)-P^{(i)}_a(b)|}
{(b-a)^{m-i}}\right\}\cdot
\,(x-a)^{m-n}
$$
and
$$
Y_2(x)=|P_b^{(n)}(x)|+
\left\{\smed_{i=0}^{m-1}\,\,
\frac{|P^{(i)}_b(a)-P^{(i)}_a(a)|}
{(b-a)^{m-i}}\right\}\cdot
\,(b-x)^{m-n}\,.
$$
\end{lemma}
\par {\it Proof.} Definition \rf{H-J} implies the existence of constants $\gamma_m,\gamma_{m+1},...,\gamma_{2m-1}\in \R$ such that  
$$
H_J(x)=P_a(x)+\,
\smed_{k=m}^{2m-1}\,\,\frac{1}{k!}\,\gamma_k\,
(x-a)^k~~~\text{for every}~~~x\in[a,b]\,.
$$
Hence, for every $n\in\{0,...,m\}$ and every $x\in[a,b]$,
\bel{DI-H}
H^{(n)}_J(x)=P^{(n)}_a(x)+\,
\smed_{k=m}^{2m-1}\,\,\frac{1}{(k-n)!}\,\gamma_k\,
(x-a)^{k-n}\,.
\ee
In particular,
$$
H^{(n)}_J(b)=P^{(n)}_a(b)+\,
\smed_{k=m}^{2m-1}\,\,\frac{1}{(k-n)!}\,\gamma_k\,
(b-a)^{k-n}
$$
which together with \rf{H-J} implies that
$$
\smed_{k=m}^{2m-1}\,\,\gamma_k\,
\frac{(b-a)^{k-n}}{(k-n)!}=P^{(n)}_b(b)-P^{(n)}_a(b),~~~
\text{for all}~~~n=0,...,m-1\,.
$$
\par Thus, the tuple $(\gamma_m,\gamma_{m+1},...,\gamma_{2m-1})$ is a solution of the above system of $m$ linear equations with respect to $m$ unknowns. Whitney \cite{W2} proved the existence of constants $K_{k,i}$,~ $k=m,...,2m-1$,~ $i=0,...,m-1$, depending only on $m$, such that
\bel{GM}
\gamma_k=\smed_{i=0}^{m-1}\,\,K_{k,i}\,
\frac{P^{(i)}_b(b)-P^{(i)}_a(b)}{(b-a)^{k-i}},~~~~~~
k=m,...,2m-1.
\ee
\par This representation enables us to estimate $H^{(n)}_J$ as follows: Thanks to \rf{GM},
$$
|\gamma_k|\le C(m)\,\smed_{i=0}^{m-1}\,\,
\frac{|P^{(i)}_b(b)-P^{(i)}_a(b)|}{(b-a)^{k-i}},~~~
k=m,...,2m-1,
$$
and, thanks to \rf{DI-H},
$$
|H^{(n)}_J(x)|\le \,|P^{(n)}_a(x)|+
\smed_{k=m}^{2m-1}\,\,\frac{1}{(k-n)!}\,|\gamma_k|\,
(x-a)^{k-n}~~~~\text{for every}~~~x\in[a,b].
$$
Hence,
\be
|H^{(n)}_J(x)|&\le& \,|P^{(n)}_a(x)|+C(m)
\smed_{k=m}^{2m-1}\,\smed_{i=0}^{m-1}\,\,
\frac{|P^{(i)}_b(b)-P^{(i)}_a(b)|}{(b-a)^{k-i}}
(x-a)^{k-n}\nn\\
&=&\,|P^{(n)}_a(x)|+C(m)
\smed_{i=0}^{m-1}\,\smed_{k=m}^{2m-1}\,\,
\frac{|P^{(i)}_b(b)-P^{(i)}_a(b)|}{(b-a)^{k-i}}
(x-a)^{k-n}\nn\\
&=&
\,|P^{(n)}_a(x)|+C(m)\smed_{i=0}^{m-1}\,
|P^{(i)}_b(b)-P^{(i)}_a(b)|\,\frac{(x-a)^{m-n}}{(b-a)^{m-i}}
\cdot\,\smed_{k=m}^{2m-1}\,\left(\frac{x-a}{b-a}\right)^{k-m}
\nn\\
&\le&
\,|P^{(n)}_a(x)|+C(m)\,m\smed_{i=0}^{m-1}\,
|P^{(i)}_b(b)-P^{(i)}_a(b)|\,\frac{(x-a)^{m-n}}{(b-a)^{m-i}}
\nn
\ee
proving that $|H^{(n)}_J(x)|\le C(m)\,m\, Y_1(x)$ for all $x\in[a,b]$. Interchanging the roles of $a$ and $b$, we show that $|H^{(n)}_J(x)|\le C(m)\, Y_2(x)$ on$[a,b]$ proving the lemma.\bx
\begin{lemma}\lbl{DF-P} Let $J=(a,b)\in\Jc_E$ be a bounded interval. Then for every $x\in [a,b]$ the following inequality
$$
|H_J^{(m)}(x)|\le C(m)\,\min\left\{\smed_{i=0}^{m-1}\,\,
\frac{|P^{(i)}_b(b)-P^{(i)}_a(b)|}{(b-a)^{m-i}},\,
\smed_{i=0}^{m-1}\,\,
\frac{|P^{(i)}_b(a)-P^{(i)}_a(a)|}{(b-a)^{m-i}}\,
\right\}
$$
holds.
\end{lemma}
\par {\it Proof.} The proof is immediate from Lemma \reff{DF-N} because $P_a$ and $P_b$ belong to $\Pc_{m-1}$.\bx
\begin{lemma}\lbl{V-J} Let $\Ic$ be a finite family of pairwise disjoint closed intervals $I=[u_I,v_I]$ such that $(u_I,v_I)\subset\R\setminus E$ for every $I\in\Ic$. Let $F$ be the function defined by \rf{DEF-F}. Then
\bel{F-L}
\smed_{I=[u_I,v_I]\in\Ic}\,
\frac{|\,F^{(m-1)}(v_I)-F^{(m-1)}(u_I)|^p}
{(v_I-u_I)^{p-1}} \le C(m)^p\, \lambda^p\,.
\ee
\end{lemma}
\par {\it Proof.} Let $I=[u_I,v_I]\in\Ic$. Since  $(u_I,v_I)\subset\R\setminus E$, there exist an interval $J=(a,b)\in\Jc_E$ containing $(u_I,v_I)$. (Recall that the family $\Jc_E$ is defined by \rf{TA-E}). The extension formula \rf{DEF-F} tells us that $F|_J=H_J$. This property and Lemma \reff{DF-P} imply that
\be
|\,F^{(m-1)}(u_I)-F^{(m-1)}(v_I)|&=&
|\,H_J^{(m-1)}(u_I)-H_J^{(m-1)}(v_I)|\le\, (\max_I |H^{(m)}|) \cdot (v_I-u_I)\nn\\
&\le&
C(m)\, (v_I-u_I)\,\smed_{i=0}^{m-1}\,\,
\frac{|P^{(i)}_a(a)-P^{(i)}_b(a)|}{(b-a)^{m-i}}\,.
\nn
\ee
Hence,
\bel{F-UV}
V_I\le
C(m)^p\,m^p\,(v_I-u_I)\,\smed_{i=0}^{m-1}\,\,
\frac{|P^{(i)}_a(a)-P^{(i)}_b(a)|^p}{(b-a)^{(m-i)p}}\,.
\ee
\par For every $J=(a,b)\in\Jc_E$ by $\Ic_J$ we denote a subfamily of $\Ic$ defined by
$$
\Ic_J=\{I\in\Ic: I\subset [a,b]\}\,.
$$
\par Let\, $\Jcw=\{J\in\Jc_E: \Ic_J\ne\emp\}$. Then, thanks to \rf{F-UV}, for every $J=(a_J,b_J)\in\Jcw$
$$
Q_J=\,\smed_{I=[u_I,v_I]\in \Ic_J} \hspace*{-2mm} V_I\,
\le
C(m)^p\,\left(\smed_{I\in \Ic_J}\diam I\right)\,\left(\smed_{i=0}^{m-1}\,\,
\frac{|P^{(i)}_{a_J}(a_J)-P^{(i)}_{b_J}(a_J)|^p}
{(b_J-a_J)^{(m-i)p}}
\right).
$$
\par We know that the intervals of the family $\Ic_J$ are pairwise disjoint (because the intervals of $\Ic$ have this property). Hence,
$$
Q_J\le
C(m)^p\,(b_J-a_J)\,\left(\smed_{i=0}^{m-1}\,\,
\frac{|P^{(i)}_{a_J}(a_J)-P^{(i)}_{b_J}(a_J)|^p}
{(b_J-a_J)^{(m-i)p}}
\right)
=
C(m)^p\,\smed_{i=0}^{m-1}\,\,
\frac{|P^{(i)}_{a_J}(a_J)-P^{(i)}_{b_J}(a_J)|^p}
{(b_J-a_J)^{(m-i)p-1}}\,.
$$
\par Finally,
$$
Q=
\smed_{I=[u_I,v_I]\in\Ic}
V_I =\smed_{J=(a_J,b_J)\in\Jcw}Q_J
\le
C(m)^p\smed_{J=(a_J,b_J)\in\Jcw}\,\,\,
\smed_{i=0}^{m-1}\,\,
\frac{|P^{(i)}_{a_J}(a_J)-P^{(i)}_{b_J}(a_J)|^p}
{(b_J-a_J)^{(m-i)p-1}}\,.
$$
\par Note that the family $\Jcw$ consist of pairwise disjoint intervals. Therefore, thanks to assumption  \rf{L-N}, $Q\le C(m)^p\,\lambda^p$, completing
the proof of the lemma.\bx
\begin{lemma}\lbl{UV-E} Let $\Ic=\{I=[u_I,v_I]\}$ be a finite family of closed intervals such that $u_I,v_I\in E$ for each $I\in\Ic$. Suppose that the open intervals $\{(u_I,v_I): I\in\Ic\}$ are pairwise disjoint. Then inequality \rf{F-L} holds.
\end{lemma}
\par {\it Proof.} Because $F$ agrees with the Whitney $(m-1)$-field $\VP=\{P_x:x\in E\}$, see \rf{F-AG-P}, we have $F^{(m-1)}(x)=P^{(m-1)}_x(x)$ for every $x\in E$. Hence,
$$
A=
\smed_{I=[u_I,v_I]\in\Ic}
\frac{|\,F^{(m-1)}(u_I)-F^{(m-1)}(v_I)|^p}
{(v_I-u_I)^{p-1}} =
\smed_{I=[u_I,v_I]\in\Ic}
\frac{|\,P^{(m-1)}_{u_I}(u_I)-P^{(m-1)}_{v_I}(v_I)|^p}
{(v_I-u_I)^{p-1}}\,.
$$
\par Because the intervals $\{(u_I,v_I): I\in\Ic\}$ are pairwise disjoint, assumption \rf{L-N} implies that $A\le \lambda^p$ proving the lemma.\bx
\medskip
\par We are in a position to finish the proof of the sufficiency. Let $\Ic$ be a finite family of pairwise disjoint closed intervals. We introduce the following notation: given an interval $I=[u,v]$, $u\ne v$, we put
$$
Y(I;F)=\frac{|\,F^{(m-1)}(u_I)-F^{(m-1)}(v_I)|^p}
{(v_I-u_I)^{p-1}}\,.
$$
We put $Y(I;F)=0$ whenever $u=v$, i.e., $I=[u,v]$ is a singleton.
\smsk
\par To each interval $I\in\Ic$ we assign three intervals $I^{(1)}$, $I^{(2)}$, $I^{(2)}$ as follows:
\smsk
\par Let $I=[u_I,v_I]\in\Ic$ be an interval such that
$I\cap E\ne\emp~~~\text{and}~~~\{u_I,v_I\}\nsubset E$.
Thus either $u_I$ or $v_I$ belongs to $\R\setminus E$.
Let $u'_I$ and $v'_I$ be the points of $E$ nearest to $u_I$ and $v_I$ on $I\cap E$ respectively. Then $[u'_I,v'_I]\subset[u_I,v_I]$. Let
$$
I^{(1)}=[u_I,u'_I],~~I^{(2)}=[u'_I,v'_I]~~~\text{and}~~~
I^{(3)}=[v'_I,v_I]\,.
$$
\par Note that $u'_I, v'_I\in E$ and $(u_I,u'_I), (v'_I,v_I)\subset \R\setminus E$ provided $u_I\notin E$ and $v_I\notin E$. Furthermore,
$$
Y(I;F)\le 3^p\,\left\{Y(I^{(1)};F)+Y(I^{(2)};F)+Y(I^{(3)};F)
\right\}.
$$
\par If $I\in\Ic$ and $(u_I,v_I)\subset\R\setminus E$, or $u_I,v_I\in E$, we put $I^{(1)}=I^{(2)}=I^{(3)}=I$.
\smsk
\par Clearly,
$$
A(F;\Ic)=\sum_{I\in\Ic}\, Y(I;F)\le
3^p\,\sum_{I\in\Ic}
\left\{Y(I^{(1)};F)+Y(I^{(2)};F)+Y(I^{(3)};F)\right\}
$$
proving that
\bel{A-IP}
A(F;\Ic)\le 3^p\,\sum_{I\in\Icw}\,Y(I;F).
\ee 
Here $\Icw=\cup\{\Icw_j:~j=1,2,3\}$  where
$\Icw_j=\cup\{I^{(j)}:~ I\in\Ic\}$,~ $j=1,2,3$.
\par We know that for each $I=[u_I,v_I]\in\Icw$ either
$(u_I,v_I)\in\R\setminus E$, or $u_I,v_I\in E$, or $u_I=v_I$ (and so $Y(I;F)=0$). Furthermore, the open intervals $\{(u_I,v_I):I\in\Icw\}$ are pairwise disjoint. This property of $\Icw$, inequality \rf{A-IP}, Lemma \reff{V-J} and Lemma \reff{UV-E} imply that  $A(F;\Ic)\le C(m)^p\,\lambda^p$.
\smsk
\par Because $\Ic$ is an arbitrary finite family of pairwise disjoint closed intervals, the function $G=F^{(m-1)}$ satisfies the hypothesis of Theorem \reff{CR-SOB}. This theorem tells us that the function $F^{(m-1)}\in \LOPR$ and $\|F^{(m-1)}\|_{\LOPR}\le C(m)\,\lambda$. This implies that $F\in\LMPR$ and $\|F\|_{\LMPR}\le C(m)\,\lambda$ proving
the sufficiency part of Theorem \reff{JET-V}.\bx
\smsk
\par  The proof of Theorem \reff{JET-V} is complete.\bx

\SECT{5. The Main Lemma: from jets to Lagrange polynomials.}{5}
\addtocontents{toc}{5. The Main Lemma: from jets to Lagrange polynomials. \hfill \thepage\par\VST}

\indent\par This section is devoted to the second main ingredient of our proof of the sufficiency part of Theorem \reff{MAIN-TH}, the Main Lemma \reff{X-SXN}. Let $E\subset\R$ be a closed set with $\#E\ge m+1$. Let $\lambda>0$ and let $f$ be a function on $E$ satisfying condition \rf{A-FE}, i.e.,
$$
\sup_{S\subset E,\,\#S=m+1}\,
|\Delta^{m}f[S]|\,(\diam S)^{\frac1p}\le \,\lambda.
$$
\par In Section 3 we have proved that in this case the Whitney field $\VP^{(m,E)}[f]=\{P_x\in\PMO:x\in E\}$ satisfying conditions \rf{FJ-H}-\rf{PX-M} of Definition \reff{P-X}, is well defined. The Main Lemma \reff{X-SXN} below provides a controlled transition from the (Hermite) polynomials $\{P_x: x\in E\}$ of the function $f$ to its Lagrange polynomials.
\begin{mlemma}\lbl{X-SXN} Let $k\in\N$, $\ve>0$, and let $X=\{x_1,...,x_k\}\subset E$, $x_1<...<x_k$, be a sequence of points in $E$. There exist a positive integer $\ell\ge m$, a finite strictly increasing sequence $V=\{v_1,...,v_\ell\}$ of points in $E$, and a mapping $H:X\to 2^V$ such that:
\msk 
\par ($\blbig 1$) For every $x\in X$ the set $H(x)$ consists of $m$ consecutive points of the sequence $V$. Thus,
$$
H(x)=\{v_{j_1(x)},...,v_{j_2(x)}\}
$$
where $1\le j_1(x)\le  j_2(x)=j_1(x)+m-1\le \ell$;
\smallskip 
\par ($\blbig 2$) $x\in H(x)$ for each $x\in X$. In particular, $X\subset V$.
\par Furthermore, given $i\in\{1,...,k-1\}$ let $x_i=v_{\vkp_i}$ and $x_{i+1}=v_{\vkp_{i+1}}$. Then $0<\vkp_{i+1}-\vkp_{i}\le 2m$.
\bigskip
\par ($\blbig 3$) Let $x',x''\in X$, $x'< x''$. Then
$$
\min H(x')\le \min H(x'')~~~~\text{and}~~~~\max H(x')\le \max H(x'')\,;
$$
\par ($\blbig 4$) For every $x', x''\in X$ such that $H(x')\ne H(x'')$ the following inequality
\bel{HX-DM1}
\diam H(x')+\diam H(x'')\le 2(m+1)\, |x'-x''|
\ee
holds;
\smallskip
\par ($\blbig 5$) For every $x,y\in X$ and every $i,0\le i\le m-1$, we have
\bel{FL-PH}
|P_{x}^{(i)}(y)-L_{H(x)}^{(i)}[f](y)|<\ve.
\ee
\end{mlemma}
\par {\it Proof.} We proceed by steps.
\medskip
\par {\bf STEP 1.} At this step we introduce the sequence $V$ and the mapping $H$.\smallskip
\par We recall that, given $x\in E$, by $\SH_x$ and $s_x$ we denote a subset of $E$ and a point in $E$ whose properties are described in Propositions \reff{SET-SX} and \reff{PR-SX}. In particular, $\#\SH_x\le m$. Let
\bel{SX-U}
\Sc_X=\bigcup_{x\in X}\, \SH_x~~~~\text{ and let}~~~~n=\#\Sc_X.
\ee
\par Clearly, one can consider $\Sc_X$ as a finite strictly increasing sequence of points $\{u_i\}_{i=1}^n$ in $E$. Thus
\bel{SXU-R}
\Sc_X=\{u_1,...,u_n\}~~~\text{and}~~~u_1<u_2<...<u_n\,.
\ee
\par If $\#\SH_x=m$ {\it for every} $x\in X$, we set $V=\Sc_X$ and $H(x)=\SH_x,$ $x\in X$. In this case the required properties ($\blbig 1$)-($\blbig 5$) of the Main Lemma are immediate from Propositions \reff{SET-SX} and \reff{PR-SX}.
\smallskip
\par However, in general, the set $X$ may have points $x$ with  $\#\SH_x<m$. For those $x$ we construct the required set $H(x)$ by adding to $\SH_x$ a certain finite set $\tH(x)\subset E$. In other words, we define $H(x)$ as
\bel{DF-HX}
H(x)=\tH(x)\cup \SH_x,~~~x\in X,
\ee
where $\tH(x)$ is a subset of $E$ such that
$\tH(x)\cap \SH_x=\emp$~ and~ $\#\tH(x)=m-\#\SH_x.$
\par Finally, we set
\bel{DF-V}
V=\cup\,\{H(x):x\in X\}\,.
\ee
\par We construct $\tH(x)$ by picking $(m-\#\SH_x)$ points of $E$ in a certain small neighborhood of $s_x$. Propositions \reff{SET-SX}, \reff{PR-SX}, and Lemma
\reff{HP-CN} enable us to prove that this neighborhood can be chosen so small that $V$ and $H$ will satisfy conditions ($\blbig 1$)-($\blbig 5$) of the Main Lemma.
\medskip
\par We turn to the precise definition of the mapping $\tH$.
\par First, we set $\tH(x)=\emp$ whenever $\#\SH_x=m$. Thus,
\bel{HW-0}
H(x)=\SH_x~~~\text{provided}~~~\# \SH_x=m\,.
\ee
\par Note, that in this case $P_{x}=L_{\SH_x}[f]=L_{H(x)}[f]$ (see \rf{PX-M}), so that \rf{FL-PH} holds vacuously.
\medskip
\par Let us define the sets $H(x)$ for all points $x\in X$ such that $\#\SH_x<m$.
\par We recall that part (iii) of Proposition \reff{PR-SX} tells us that
\bel{SX-LP}
\text{for each\, $x\in X$\, with\, $\#\SH_x<m$\, the point \,$s_x$\, is a {\it limit  point} of\, $E$}.
\ee
In turn, part (i) of this proposition tells us that
$$
\text{either}~~s_x=\min \SH_x~~\text{or}~~s_x=\max \SH_x\,.
$$
\par Let
$$
\LIM_X=\{s_x:x\in X,~ \#\SH_x<m\}.
$$
Then, thanks to \rf{SX-LP},
\bel{SZ-A}
\text{every point}~~~z\in \LIM_X ~~~\text{is a {\it limit  point} of}~~E.
\ee
\par Given $z\in \LIM_X$, let
\bel{KZ-D}
\KZ=\{x\in X: s_x=z\}\,.
\ee
\par The following lemma describes main properties of the sets $\KZ$, $z\in \LIM_X$.
\begin{lemma}\lbl{KZ-PR}  Let $z\in \LIM_X$. Suppose that $\KZ\ne\{z\}$. Then the following properties hold:
\smsk
\par ($\squ 1$) The set $\KZ$ lies on one side of $z$, i.e.,
\bel{SIDE-X}
\text{either}~~~\max \KZ\le z
~~~\text{or}~~~\min \KZ\ge z\,;
\ee
\par ($\squ 2$) If $\min \KZ\ge z$\,, then for every $r>0$ the interval
\bel{Z-MIN}
(z-r,z)~~~\text{contains an infinite number of points of}~~E\,.
\ee
\par If $\max \KZ\le z$ then each interval $(z,z+r)$ contains an infinite number of points of $E$;
\msk
\par ($\squ 3$) If $\min \KZ\ge z$ then
\bel{KY-S}
[z,y]\cap X\subset \KZ~~~\text{for every}~~~y\in \KZ\,.
\ee
Furthermore,
\bel{Z-MSX}
\min\SH_y=z~~\text{for all}~~y\in\KZ,
\ee
and
\bel{SX-Y}
\SH_{x}\subset \SH_{y}~~\text{for every}~~y\in \KZ~~\text{and every}~~x\in[z,y]\cap E.
\ee
\par Correspondingly, if\, $\max \KZ\le z$\, then\, $[y,z]\cap X\subset \KZ$\, for each $y\in \KZ$. Moreover, $\SH_{x}\subset \SH_{y}$ for every $y\in \KZ$ and every $x\in[y,z]\cap E$. In addition, $\max\SH_y=z$~ for all $y\in\KZ$;
\msk
\par ($\squ 4$) Assume that $\min \KZ\ge z$. Let $\brz=\max \KZ$. Then
\bel{KZ-IX}
\KZ=[z,\brz]\cap X.
\ee
\par Furthermore, in this case $\KZ\subset S_{\brz}$.
\par In turn, if $\max \KZ\le z$ then $\KZ=[\tz,z]\cap X$ where $\tz=\min \KZ$. In this case $\KZ\subset S_{\tz}$\,;
\msk
\par ($\squ 5$) $\#\KZ\le m$.
\end{lemma}
\par {\it Proof.} ($\squ 1$) Suppose that \rf{SIDE-X} does not hold so that there exist $z',z''\in \KZ$ such that $z''<z<z'$.
\par Thanks to \rf{KZ-D}, $z=s_{z''}=s_{z'}$. We also know
that $z',s_{z'}\in\SH_{z'}$, see part (i) of Proposition \reff{SET-SX} and part (i) of Proposition \reff{PR-SX}. This property and \rf{MNM-X1} tell us that
\bel{Z-SHZ}
[z,z']\cap E=[s_{z'},z']\cap E\subset [\min \SH_{z'},\max \SH_{z'}]\cap E=\SH_{z'}.
\ee
\par Part (i) of Proposition \reff{SET-SX} also tells us that $\#\SH_{z'}\le m$ proving that the interval
\bel{Z-ME}
(z,z')~~~\text{contains at most $m$ points of $E$}.
\ee
In the same way we show that the interval $(z'',z)$ contains at most $m$ points of $E$.
\par Thus, the interval $(z'',z')$ contains a finite number of points of $E$ proving that $z$ is an {\it isolated point} of $E$. On the other hand, $z\in \LIM_X$ so that, thanks to \rf{SZ-A}, $z$ is a {\it limit point} of $E$, a contradiction.
\smallskip
\par ($\squ 2$) Let $z'\in \KZ$, $z'\ne z$. Then $z<z'$ so that, thanks to \rf{Z-ME}, the interval $(z,z')$ contains at most $m$ points of $E$. But $z$ is a limit point of $E$, see  \rf{SZ-A}, so that the interval $(z-r,z)$ contains an infinite number of points of $E$. This proves \rf{Z-MIN}.
\par In the same fashion we prove the second statement of part ($\squ 2$).
\smallskip
\par ($\squ 3$) Let $\min \KZ\ge z$, and let $y\in \KZ$, $y\ne z$. Prove that $s_x=z$ for every $x\in [z,y]\cap E$.
\par We know that $z=s_y<y$. Furthermore, property \rf{Z-SHZ} tells us that
\bel{Z-YN}
[z,y]\cap E\subset \SH_{y}.
\ee
\par We recall that, thanks to \rf{SZ-A}, $z$ is a limit point of $E$ so that $\SH_z=\{z\}$. Part (iii) of Proposition \reff{SET-SX} tells us that
$z=\min\SH_z\le\min \SH_x$~ and $\max\SH_x\le\max \SH_y$,
so that
\bel{SX-Y1}
\SH_{x}\subset [z,\max \SH_{y}]\cap E.
\ee
In particular, $z\le\min \SH_y$. But $z=s_y\in \SH_{y}$, so that $z=\min \SH_{y}$ proving \rf{Z-MSX}.
\par In turn, thanks to \rf{MNM-X1},
$$
[\min \SH_{y},\max \SH_{y}]\cap E=[z,\max \SH_{y}]\cap E=\SH_{y}.
$$
This and \rf{SX-Y1} imply that $\SH_{x}\subset \SH_{y}$ for every $x\in[z,y]\cap E$ proving \rf{SX-Y}.
\par Moreover, thanks to \rf{SX-LP}, if $\#\SH_{x}<m$ then $s_x$ is a limit point of $E$. But $s_x\in\SH_x\subset \SH_{y}$, therefore, part (ii) of Proposition \reff{PR-SX} implies that $s_x=s_{y}=z$.
\smsk
\par If $\#\SH_{x}=m$ then $\SH_{x}=\SH_{y}$ (because $\SH_{x}\subset \SH_{y}$ and $\#\SH_{y}\le m$). Hence, $z=s_y\in \SH_x$. But $z$ is a limit point of $E$ which together with part (ii) of Proposition \reff{PR-SX} implies that $s_x=z=s_y$.
\par Thus, in all cases $s_x=z$ proving property ($\squ 3$) of the lemma in the case under consideration. Using the same ideas we prove the second statement of the lemma related to the case $\max \KZ\le z$.
\smallskip
\par ($\squ 4$) From \rf{KY-S}, we have
$[z,\brz]\cap X\subset \KZ$.

On the other hand, $\KZ\subset [z,\brz]$ because
$z\le \min \KZ$ and $\brz=\max \KZ$. Since $\KZ\subset X$,
see \rf{KZ-D}, $\KZ\subset [z,\brz]\cap X$ proving \rf{KZ-IX}.
\par Furthermore, thanks to \rf{Z-YN} (with $y=\brz$), it follows that $[z,\brz]\cap E\subset \SH_{\brz}$ proving that
$$
\KZ=[z,\brz]\cap X\subset [z,\brz]\cap E\subset \SH_{\brz}\,.
$$
\par In the same way we prove the last statement of part ($\squ 4$) related to the case $z\ge \max \KZ$.
\smallskip
\par ($\squ 5$) We recall that $\#\SH_{x}\le m$ for every $x\in E$, see part (i) of Proposition \reff{SET-SX}. Part ($\squ 4$) of the present lemma tells us that $\KZ\subset \SH_{y}$ where $y=\max\KZ$ or $y=\min\KZ$. Hence $\#\KZ\le \#\SH_{y}\le m$.
\smsk
\par The proof of the lemma is complete.\bx
\msk
\par Let us fix a point $z\in \LIM_X$ and define the set $H(x)$ for every $x\in \KZ$. Thanks to property \rf{SIDE-X}, it suffices to consider the following three cases:\medskip
\par {\it Case} ($\bigstar$1).~ Suppose that
\bel{C-BS1}
\KZ\ne\{z\}~~~\text{and}~~~\min \KZ\ge z.
\ee
\par Let $y\in \KZ$. Thus $y\in X$ and $s_y=z$; we also know that $y\ge z$. Then property \rf{KY-S} tells us that
$[z,y]\cap X\subset \KZ$.
\smallskip
\par We also note that $\KZ=[z,\brz]\cap X$ where  $\brz=\max \KZ$, and $\KZ\subset S_{\brz}$, see part ($\squ 4$) of Lemma \reff{KZ-PR}. Furthermore, part ($\squ 5$) of this lemma tells us that  $\#\KZ\le m$.
\medskip
\par Let us fix several positive constants which we need for definition of the sets $\{H(x):x\in X\}$.
\par We recall that $X=\{x_1,...x_k\}$ and  $x_1<...<x_k$.
Let $I_X=[x_1,x_k]$.
\par We also recall that inequality \rf{A-FE} holds, and $s_x=z$ provided $x\in \KZ$. This enables us to apply Lemma \reff{HP-CN} to the interval $I=I_X$ and the point $x\in \KZ$. This lemma tells us that
$$
\lim_{\substack{S'\setminus \SH_x\to z \smallskip\\
\SH_x\subset\, S'\subset E,\,\,\#S'=m}} \,\,\|L_{S'}[f]-P_x\|_{C^m(I_X)}=0\,.
$$
Thus, there exists a constant $\dw_x=\dw_x(\ve)>0$ satisfying the following condition: for every $m$-point set $S'$ such that $\SH_x\subset S'\subset E$ and $S'\setminus \SH_x\subset(z-\dw_x,z+\dw_x)$ we have
$$
|P_x^{(i)}(y)-L_{S'}^{(i)}[f](y)|<\ve~~~\text{for every}~~~i=0,...,m-1,~~~\text{and every}~~~y\in I_X\,.
$$
\par We recall that $\Sc_X=\{u_1,...,u_n\}$ is the set defined by \rf{SX-U} and \rf{SXU-R}.  Let
\bel{TAU-X}
\tau_X=\tfrac14\, \min_{i=1,...,n-1}(u_{i+1}-u_{i}).
\ee
Thus,
\bel{X-TAU}
|x-y|\ge 4\tau_X~~~\text{provided}~~~x,y\in \Sc_X, x\ne y\,.
\ee
\par Finally, we set
\bel{DL-Z}
\delta_z=\min\,\{\tau_X,\min_{x\in \KZ}\dw_x\}\,.
\ee
Clearly, $\delta_z>0$ (because $\KZ$ is finite).
\smallskip
\par Definition \rf{DL-Z} implies the following: Let $x\in \KZ$. Then for every $i,~0\le i\le m-1$, and every $m$-point set $S'$ such that
$\SH_x\subset S'\subset E$~ and $S'\setminus \SH_x\subset(z-\delta_z,z+\delta_z)$, we have
\bel{LM-2}
|P_x^{(i)}(y)-L_{S'}^{(i)}[f](y)|<\ve,~~~y\in X.
\ee
\par Inequality \rf{Z-MIN} tells us that the interval $(z-\delta_z,z)$ contains an infinite number of points of $E$. Let us pick $m-1$ distinct points $a_1<a_2<...<a_{m-1}$ in  $(z-\delta_z,z)\cap E$ and set
$$
W(z)=\{a_1,a_2,...,a_{m-1}\}.
$$
Thus,
\bel{W-ZD}
W(z)=\{a_1,a_2,...,a_{m-1}\}\subset (z-\delta_z,z)\cap E.
\ee
\par In particular,
\bel{AC-1}
z-\tau_X<a_1<a_2<...<a_{m-1}<z~~~~~~~(\text{because}~~~ \delta_z\le\tau_X).
\ee
\par Let $x\in \KZ$, and let $\ell_x=\# \SH_x$. We introduce the set $\tH(x)$ as follows: we set
\bel{HW-M}
\tH(x)=\emp~~~\text{and}~~~H(x)=\SH_x~~~
\text{provided}~~~\ell_x=m.
\ee
\par
If $\ell_x<m$, we define $\tH(x)$ by letting
\bel{HW-1}
\tH(x)=\{a_{\ell_x},a_{\ell_x+1},...,a_{m-1}\}.
\ee
Clearly, $\#\tH(x)+\#\SH_x=m$.
\par Then we define $H(x)$ by formula \rf{DF-HX}, i.e., we set $H(x)=\tH(x)\cup \SH_x$. This definition, property \rf{W-ZD} and inequality \rf{LM-2} imply that for every $y\in X$ and every $i,~0\le i\le m-1$, the following inequality
\bel{HX-C1}
|P_{x}^{(i)}(y)-L_{H(x)}^{(i)}[f](y)|<\ve
\ee
holds. Furthermore, property \rf{SX-Y} tells us that
$\SH_{x'}\subset \SH_{x''}$ provided $x',x''\in \KZ$, $x'< x''$. This property and definition \rf{HW-1} imply that 
\bel{C1-P3}
\min H(x')\le \min H(x'')~~~\text{for every}~~~x',x''\in \KZ,~x'< x''.
\ee
\par Let us also note the following property of $H(x)$ which directly follows from its definition: Let
\bel{H-SX}
\HH(x)=\tH(x)\cup \{s_x\}\,.
\ee
Then
\bel{H-GS}
[\min H(x),\max H(x)]=[\min\SH_x,\max\SH_x]\cup[\min \HH(x),\max \HH(x)]
\ee
and
\bel{H-MH2}
[\min\SH_x,\max\SH_x]\cap[\min \HH(x),\max \HH(x)]=\{s_x\}.
\ee
\smallskip
\par  {\it Case} ($\bigstar$2). $\KZ\ne\{z\}$~ and $\max \KZ\le z$.
\smsk
\par Using the same approach as in {\it Case} ($\bigstar$1), see \rf{C-BS1}, given $x\in \KZ$, we define a corresponding constant $\delta_z$, a set $W(z)=\{a_1,...,a_{m-1}\}$ and sets $\tH(x)$ and $H(x)$. More specifically, we pick  a strictly increasing sequence
$$
W(z)=\{a_1,a_2,...,a_{m-1}\}\subset (z,z+\delta_z)\cap E.
$$
In particular, this sequence has the following property:
\bel{AC-2}
z<a_1<a_2<...<a_{m-1}<z+\tau_X.
\ee
(Recall that $\tau_X$ is defined by \rf{TAU-X}.) Then we set $\tH(x)=\emp$ and $H(x)=\SH_x$ if $\ell_x=m$, and
\bel{HW-2}
\tH(x)=\{a_1,a_2,...,a_{m-\ell_x}\}~~~\text{if}~~~
\ell_x<m.
\ee
(Recall that $\ell_x=\#\SH_x$.) Finally, we define the set $H(x)$ by formula \rf{DF-HX}.
\par As in {\it Case} ($\bigstar$1), see \rf{C-BS1}, our choice of $\delta_z$, $\tH(x)$ and $H(x)$ provides inequality \rf{HX-C1} and pro\-perties \rf {C1-P3}, \rf{H-GS} and \rf{H-MH2}.
\bigskip
\par  {\it Case} ($\bigstar$3). $\KZ=\{z\}$.
\smsk
\par Note that in this case $z\in X$ is a limit point of $E$ and $\SH_z=\{z\}$. This enables us to pick a subset
$W(z)=\{a_1,...,a_{m-1}\}\subset E$ with~ $\#W(z)=m-1$ such that either \rf{AC-1} or \rf{AC-2} hold.
\par We set $\tH(z)=W(z)$. Thus, in this case the set $H(z)$ is defined by formula \rf{DF-HX} with
\bel{HW-3}
\tH(z)=\{a_1,a_2,...,a_{m-1}\},
\ee
i.e., $H(z)=\{z,a_1,a_2,...,a_{m-1}\}$.
\par It is also clear that properties \rf{H-GS}, \rf{H-MH2} hold in the case under consideration. Moreover, our choice of the set $W(z)$ provides  inequality \rf{HX-C1} with $x=z$ and $H(x)=\{z,a_1,a_2,...,a_{m-1}\}$.
\bigskip
\par We have defined the set $H(x)$ for every $x\in X$. Then we define the set $V$ by formula \rf{DF-V}. Clearly, $V$ is a finite subset of $E$. Let us enumerate the points of this set in increasing order: thus, we represent $V$ in the form
$$
V=\{v_1,v_2,...,v_\ell\}
$$
where $\ell$ is a positive integer and $\{v_j\}_{j=1}^\ell$ is a strictly increasing sequence of points in $E$.
\bigskip
\par {\bf STEP 2.} At this step we prove two auxiliary lemmas which describe a series of important properties of the mappings $\tH$ and $H$.
\begin{lemma}\lbl{H-U1} (i) For each $x\in X$ the following inclusion
\bel{HW-IS}
\HH(x)\subset(s_x-\tau_X,s_x+\tau_X)
\ee
holds. (Recall that $\HH(x)=\tH(x)\cup \{s_x\}$, see \rf{H-SX}.)
\par (ii) The following property
$$
[\min H(x),\max H(x)]\subset [\min \SH(x)-\tau_X,\max \SH(x)+\tau_X]
$$
holds for every $x\in X$.
\end{lemma}
\par {\it Proof.} Property (i) is immediate from \rf{AC-1}, \rf{H-SX}, \rf{AC-2}. In turn, property (ii) is immediate from \rf{SX-S}, \rf{DF-HX} and \rf{HW-IS}.\bx
\begin{lemma}\lbl{H-U2} Let $x,y\in X$. Suppose that  $\#\SH_x<m$ and
\bel{TH-S}
[\min\HH(x),\max\HH(x)]\cap
[\min H(y),\max H(y)]\ne\emp\,.
\ee
Then $s_x=s_y$.
\end{lemma}
\par {\it Proof.} Part (iii) of Proposition \reff{PR-SX} tells us that the point $s_x$ {\it is a limit point of} $E$. Thus, if $y=s_x$, part (ii) of Proposition \reff{PR-SX} implies that $s_y=y$ so that in this case the lemma holds.
\smsk
\par Let us prove the lemma for $y\ne s_x$. To do so, we assume that $s_x\ne s_y$.
\par First, we prove that
$s_x\notin [\min \SH_y,\max \SH_y]$. Indeed, suppose that $s_x\in [\min \SH_y,\max \SH_y]$. In this case part (i) of Proposition \reff{SET-SX} tells us that
$$
s_x\in [\min \SH_y,\max \SH_y]\cap E=\SH_y.
$$
But $s_x$ is a limit point of $E$, so that, thanks to part (ii) of Proposition \reff{PR-SX}, $s_x=s_y$ which  contradicts our assumption $s_s\ne s_y$.
\par Thus, $s_x\notin [\min \SH_y,\max \SH_y]$ so that  $s_x\ne \min \SH_y$ and $s_x\ne \max \SH_y$. On the other hand, thanks to \rf{SX-U}, the points  $s_x,\min \SH_y,\max \SH_y$ belong to $\Sc_X$. Therefore, thanks to \rf{X-TAU},
$$
|s_x-\min\SH_y|\ge 4\tau_X~~~~~\text{and}~~~~~|s_x-\max\SH_y|\ge 4\tau_X.
$$
Hence,
$$
\dist(s_x,[\min \SH_y,\max \SH_y])\ge 4\tau_X.
$$
\par On the other hand, part (i) and part (ii) of Lemma \reff{H-U1} tell us that $\HH(x)\subset(s_x-\tau_X,s_x+\tau_X)$
and
$$
[\min H(y),\max H(y)]\subset
[\min \SH_y-\tau_X,\max \SH_y+\tau_X].
$$
Hence,
$$
[\min\HH(x),\max\HH(x)]\cap
[\min H(y),\max H(y)]=\emp.
$$
This contradicts \rf{TH-S} proving that the assumption $s_x\ne s_y$ does not hold.\bx
\bigskip
\par {\bf STEP 3.} We are in a position to prove properties ($\blbig 1$)-($\blbig 5$) of the Main Lemma \reff{X-SXN}.
\medskip
\par $\blacksquare$~ {\it Proof of property ($\blbig 1$).}  This property is equivalent to the following statement:
\bel{H-V1}
[\min H(x),\max H(x)]\cap V=H(x)~~~~~\text{for every} ~~~~x\in X.
\ee
\par Let us assume that \rf{H-V1} does not hold for certain $x\in X$, and show that this assumption leads to a contradiction.
\par Thanks to definition \rf{DF-V}, if \rf{H-V1} does not hold  then there exist $y\in X$ and $u\in H(y)$ such that
\bel{U-HX}
u\in [\min H(x),\max H(x)]\setminus H(x).
\ee
\par Prove that $\#\SH_x<m$. Indeed, otherwise, $\SH_x=H(x)$ (see \rf{HW-0}). In this case \rf{MNM-X1} implies that $[\min H(x),\max H(x)]\cap E=H(x)$ so that
$[\min H(x),\max H(x)]\setminus H(x)=\emp$. This contradicts \rf{U-HX} proving that $\#\SH_x<m$.
\par Furthermore, property \rf{MNM-X1} tells us that
$[\min \SH_x,\max\SH_x]\cap E=\SH_x\subset H(x)$
so that
\bel{U-3}
u\notin [\min \SH_x,\max\SH_x].
\ee
From this and \rf{H-GS}, we have
\bel{U-4}
u\in [\min\HH(x),\max\HH(x)].
\ee
\par We conclude that the hypothesis of Lemma \reff{H-U2} holds for $x$ and $y$ because $\#\SH_x<m$, $u\in H(y)$ and \rf{U-4} holds. This lemma tells us that $s_x=s_y$.
\smsk
\par Let $z=s_x$ so that $x,y\in \KZ$, see \rf{KZ-D}. Clearly, $\KZ\ne\{z\}$; otherwise $x=y$ contradicting \rf{U-HX}.
\par We know that either $\min \KZ\ge z$ (i.e., $z$ satisfies the condition of {\it Case} ($\bigstar$1) of\, {\bf STEP 1}, see \rf{C-BS1}), or $\max \KZ\le z$~ ({\it Case} ($\bigstar$2) of\, {\bf STEP 1} holds).
\smsk
\par Suppose that $\min \KZ\ge z$, see \rf{C-BS1}. Then $\min \SH_x=s_x=z$ so that
\bel{U-11}
[\min \SH_x,\max\SH_x]=[z,\max\SH_x].
\ee
Moreover, thanks to \rf{HW-1} and \rf{H-SX},
\bel{HTH}
\tH(x)=\{a_{\ell_x},a_{\ell_x+1},...,a_{m-1}\} ~~~\text{and}~~~
\HH(x)=\{a_{\ell_x},...,a_{m-1},z\}.
\ee
Here $\ell_x=\# \SH_x$ and $a_{1},...,a_{m-1}$ are $m-1$ distinct points of $E$ satisfying inequality \rf{AC-1}.
\par Properties \rf{U-3} and \rf{U-11} tell us that
$u\ne z=s_x$. From this, \rf{U-4} and \rf{HTH}, we have $u\in [a_{\ell_x},z)$. On the other hand, $u\in H(y)$. Because $y\in \KZ$, definitions \rf{DF-HX} and \rf{HW-1} imply that
$$
H(y)=\{a_{\ell_y},...,a_{m-1}\}\cup \SH_y.
$$
\par But $\min \SH_y=z$, see \rf{Z-MSX}. We also know that $u\in [a_{\ell_x},z)$, so that  $u\in\{a_{\ell_y},...,a_{m-1}\}$ and $u\in [a_{\ell_x},z)$. Hence, $u\in\{a_{\ell_x},...,a_{m-1}\}\subset H(x)$ which contradicts \rf{U-HX}.
\par In the same way we obtain a contradiction if $z$ satisfies the condition of {\it Case} ($\bigstar$2) of\, {\bf STEP 1}.
\smsk
\par The proof of property ($\blbig 1$) of the Main Lemma is complete.
\msk
\par $\blacksquare$~ {\it Proof of property ($\blbig 2$).} Part (i) of Proposition \reff{PR-SX} tells us that $x\in \SH_x$ for every $x\in E$. In turn, definition \rf{DF-HX} implies that $\SH_x\subset H(x)$ so that $x\in H(x)$. Hence, $x\in V$ for each $x\in X$, see \rf{DF-V}, proving that $X\subset V$.
\par Let us prove that $0<\vkp_{i+1}-\vkp_{i}\le 2m$ provided $x_i=v_{\vkp_i}$ and $x_{i+1}=v_{\vkp_{i+1}}$. The first inequality is obvious because $x_i<x_{i+1}$ and $V=\{v_j\}_{j=1}^\ell$ is a strictly increasing sequence.
\par Our proof of the second inequality relies on the following fact:
\bel{V-HH}
V\cap[x_i,x_{i+1}]\subset H(x_i)\cup H(x_{i+1}).
\ee
\par Indeed, let $v\in V\cap[x_i,x_{i+1}]$. Then definition \rf{DF-V} implies the existence of a point $\xb\in X$ such that $H(\xb)\ni v$. Hence, $v\le \max H(\xb)$.
\par Suppose that $\xb<x_i$. In this case property
($\blbig 3$) of the Main Lemma \reff{X-SXN} (which we prove below) tells us that $\max H(\xb)\le \max H(x_i)$ so that $x_i\le v\le \max H(x_i)$. Hence, $v\in [\min H(x_i),\max H(x_i)]$ (because $x_i\in H(x_i)$). We also know that $v\in V$. This and property \rf{H-V1} (which is equivalent to property ($\blbig 1$) of the Main Lemma) imply that
$$
v\in [\min H(x_i),\max H(x_i)]\cap V=H(x_i).
$$
\par In the same way we show that $v\in H(x_{i+1})$ provided $\xb>x_{i+1}$, and the proof of \rf{V-HH} is complete.
\smsk
\par Because $\#H(x_i)=\#H(x_{i+1})=m$, property \rf{V-HH} tells us that the interval $[x_i,x_{i+1}]$ contains at most $2m$ points of the set $V$. This implies the required second inequality $\vkp_{i+1}-\vkp_{i}\le 2m$ completing the proof of part ($\blbig 2$) of the Main Lemma.
\medskip
\par $\blacksquare$~ {\it Proof of property ($\blbig 3$).} Let $x',x''\in X$, $x'<x''$. Prove that
\bel{HXX}
\min H(x')\le\min H(x'').
\ee
\par We recall that
$$
H(x')=\tH(x')\cup \SH_{x'}~~~\text{and}~~~H(x'')=\tH(x'')\cup \SH_{x''},
$$
see \rf{DF-HX}. Here $\tH$ is the set defined by formulae \rf{HW-M}, \rf{HW-1}, \rf{HW-2} and \rf{HW-3}.
\par Part (iii) of Proposition \reff{SET-SX} tells us that $\min \SH_{x'}\le\min \SH_{x''}$. Since $s_{x''}\in\SH_{x''}$, we have
\bel{A-6}
\min \SH_{x'}\le\min \SH_{x''}\le s_{x''}\,.
\ee
\par We proceed the proof of \rf{HXX} by cases.
\smsk
\par {\it Case 1.} Assume that $\min \SH_{x'}< s_{x''}$.
\par Because the points $\min \SH_{x'}$ and $s_{x''}$ belong to the set $\Sc_X$ (see \rf{SX-U}), inequality \rf{X-TAU} implies that
\bel{A-8}
s_{x''}-\SH_{x'}>4\tau_X\,.
\ee
On the other hand, part (i) of Lemma \reff{H-U1} tells us that
$$
\tH(x'')\subset \HH(x'')\subset (s_{x''}-\tau_X,s_{x''}+\tau_X).
$$
(Recall that $\HH(x'')=\tH(x'')\cup \{s_{x''}\}$, see \rf{H-SX}.) From this inclusion and \rf{A-8}, we have
$$
\min\tH(x'')>s_{x''}-\tau_X>\min \SH_{x'}+3\tau_X>\min \SH_{x'}\,.
$$
This and \rf{A-6} implies us that
\be
\min H(x'')&=&\min\,(\tH(x'')\cup \SH_{x''})=
\min\,\{\min\tH(x''),\,\min \SH_{x''}\}\nn\\
&\ge&\min \SH_{x'}
\ge \min\,(\tH(x')\cup\SH_{x'})=\min H(x')\nn
\ee
proving \rf{HXX} in the case under consideration.
\msk
\par {\it Case 2.} Suppose that
$\min \SH_{x'}=s_{x''}$, and consider two cases.
\smsk
\par {\it Case 2.1:}\, $\#\SH_{x''}=m$. Then $H(x'')=\SH_{x''}$, see \rf{HW-0}, which together with  \rf{A-6} implies that
$$
\min H(x')\le \min \SH_{x'}\le \min \SH_{x''}=\min H(x'')
$$
proving \rf{HXX}.
\smsk
\par {\it Case 2.2:}\, $\#\SH_{x''}<m$. In this case part (iii) of Proposition \reff{PR-SX} tells us that $s_{x''}$ is a limit point of $E$. We know that $s_{x''}=\min \SH_{x'}$ so that the point $\min \SH_{x'}$ is a limit point of $E$ as well. Hence, $\min \SH_{x'}=s_{x'}$, see part (ii) of Proposition \reff{PR-SX}.
\smsk
\par Thus, $s_{x'}=s_{x''}=\min \SH_{x'}$. Let $z=s_{x'}=s_{x''}$. We know that
$z=\min \SH_{x'}\le x'<x''$, so that $x',x''\in \KZ$, see \rf{KZ-D}. In particular, $\KZ\ne\{z\}$ proving that the point $z$ satisfies the condition of {\it Case} ($\bigstar$1) of {\bf STEP 1}, see \rf{C-BS1}. In this case inequality \rf{HW-1} tells us that $\min H(x')\le \min H(x'')$ proving \rf{HXX} in {\it Case 2.2}.
\smallskip
\par Thus, \rf{HXX} holds in all cases. Since each of the sets $H(x')$ and $H(x'')$ consists of $m+1$ consecutive points of the strictly increasing sequence $V=\{v_i\}_{i=1}^\ell$ and $\min H(x')\le\min H(x'')$, we conclude that  $\max H(x')\le\max H(x'')$. 
\par The proof of property ($\blbig 3$) is complete.\medskip
\par $\blacksquare$~ {\it Proof of property ($\blbig 4$).} Let $x',x''\in X$, $x'\ne x''$, and let
$H(x')\ne H(x'')$.
\par Part (i) of Lemma \reff{SET-SX} and definition \rf{SX-U} tell us that
$x'\in \SH_{x'}, x''\in \SH_{x''}$ and $\SH_{x'}, \SH_{x''}\in\Sc_X$.
Hence, $x',x''\in\Sc_X$ so that, thanks to
\rf{X-TAU}, $|x'-x''|\ge 4\tau_X$. In turn, definition \rf{DF-HX} implies that
$$
\diam  H(x')\le\diam \tH(x')+\diam \SH_{x'}~~~\text{and}~~~ \diam  H(x'')\le\diam \tH(x'')+\diam \SH_{x''}\,.
$$
\par Furthermore, part (i) of Lemma \reff{H-U1} tells us that
$$
\max\,\{\diam\tH(x'),\diam\tH(x'')\}\le 2\tau_X<|x'-x''|.
$$
Hence,
\bel{DM-H1}
\diam H(x')\le \diam \SH_{x'}+|x'-x''|~~~\text{and}~~~
\diam H(x'')\le \diam \SH_{x''}+|x'-x''|\,.
\ee
\par Prove that $\SH_{x'}\ne \SH_{x''}$. Indeed, suppose that $\SH_{x'}=\SH_{x''}$ and prove that this equality contradicts to the assumption that $H(x')\ne H(x'')$.
\smsk
\par If $\#\SH_{x'}=\#\SH_{x''}=m$ then $\SH_{x'}=H_{x'}$ and $\SH_{x''}=H_{x''}$, see \rf{HW-0}, which implies the required contradiction $H(x')=H(x'')$.
\par Let now $\#\SH_{x'}=\#\SH_{x''}<m$. In this case $s_{x'}$ is the {\it unique} limit point of $E$ which belongs to $\SH_{x'}$, see part (ii) of Proposition \reff{PR-SX}. A similar statement is true for $s_{x''}$ and $\SH_{x''}$. Hence $s_{x'}=s_{x''}$.
\smallskip
\par Let $z=s_{x'}=s_{x''}$. Thus, $x',x''\in \KZ$, see \rf{KZ-D}. Since $x'\ne x''$, we have $\KZ\ne\{z\}$, so that either the point $z$ satisfies the condition of
{\it Case} ($\bigstar$1) of {\bf STEP 1} (see \rf{C-BS1}) or the condition of {\it Case} ($\bigstar$2) of {\bf STEP 1} holds. Furthermore, since $\#\SH_{x'}=\#\SH_{x''}<m$,
the sets $\tH(x'),\tH(x'')$ are determined by the formula \rf{HW-1} or \rf{HW-2} respectively. In both cases the definitions of the sets $\tH(x'),\tH(x'')$ depend only on the point $z$ (which is the same for $x'$ and $x''$ because $z=s_{x'}=s_{x''}$) and the number of points in the sets $\SH_{x'}$ and $\SH_{x''}$ (which of course is also the same because $\SH_{x'}=\SH_{x''}$).
\par Thus, $\tH(x')=\tH(x'')$ proving that
$$
H(x')=\tH(x')\cup \SH_{x'}=\tH(x'')\cup \SH_{x''}=H(x''),
$$
a contradiction. This contradiction proves that $\SH_{x'}\ne \SH_{x''}$. In this case part (ii) of Proposition \reff{SET-SX} tells us that
$$
\diam \SH_{x'}+\diam \SH_{x''}\le 2\,m\,|x'-x''|\,.
$$
Combining this inequality with \rf{DM-H1}, we obtain  the required inequality \rf{HX-DM1} proving the property ($\blbig 4$) of the Main Lemma.
\medskip
\par $\blacksquare$~ {\it Proof of property
($\blbig 5$).} In the process of constructing of the sets $H(x), x\in X$, we have noted that in all cases of {\bf STEP 1} ({\it Case} ($\bigstar$1) (see \rf{C-BS1}), {\it Case} ($\bigstar$2), {\it Case} ($\bigstar$3)) inequality \rf{HX-C1} holds for all $y\in X$ and all $i,~0\le i\le m-1$. This inequality coincides with inequality \rf{FL-PH} proving property ($\blbig 5$) of the Main Lemma.
\smsk
\par The proof of Main Lemma \reff{X-SXN} is complete.\bx

\SECT{6. The variational extension criterion: sufficiency.}{6}
\addtocontents{toc}{6. The variational extension criterion: sufficiency. \hfill\thepage\par \VST}

\indent\par In this section we prove the sufficiency part of Theorem \reff{MAIN-TH}. Let $E$ be a closed subset of $\R$ with $\#E\ge m+1$, and let $f$ be a function on $E$ such that $\lambda=\NMP(f:E)<\infty$. See \rf{NR-TR}. This enables us to make the following
\begin{assumption}\lbl{A-1} For every integer $n\ge m$ and every strictly increasing sequence of points $\{x_0,...,x_n\}$ in $E$, the following inequality
\bel{N-LMR}
\smed_{i=0}^{n-m}
\,\,
(x_{i+m}-x_i)\,|\Delta^mf[x_i,...,x_{i+m}]|^p
\le \lambda^p
\ee
holds.
\end{assumption}
\par Our aim is to prove that there exists a function $F\in\LMPR$ such that $F|_E=f$ and $\|F\|_{\LMPR}\le C(m)\,\lambda$.
\par Clearly, thanks to \rf{N-LMR},
$$
\sup_{S\subset E,\,\#S=m+1}\,
|\Delta^{m}f[S]|\,(\diam S)^{\frac1p}\le \,\lambda,
$$
proving that inequality \rf{A-FE} holds. As we have shown in Section 3, in this case the Whitney field
$\VP^{(m,E)}[f]=\{P_x\in\PMO:x\in E\}$
introduced in Definition \reff{P-X}, is well defined.
\smallskip
\par We prove the existence of the function $F$ with the help of Theorem \reff{JET-V} which we apply to the field $\VP^{(m,E)}[f]$. To enable us to do this, we first have to check that the hypothesis of this theorem holds, i.e., we must show that for every integer $k>1$ and every strictly increasing sequence $\{x_j\}_{j=1}^k$ in $E$ the following inequality
\bel{VP-LM}
\smed_{j=1}^{k-1}\,\,
\smed_{i=0}^{m-1}\,\,
\frac{|P^{(i)}_{x_{j+1}} (x_j)-P^{(i)}_{x_j}(x_j)|^p}
{(x_{j+1}-x_{j})^{(m-i)p-1}}\le\,C(m)^p\,\lambda^p
\ee
holds.
\medskip
\par Our proof of inequality \rf{VP-LM} relies on
Main Lemma \reff{X-SXN}. More specifically, we fix $\ve>0$ and apply this lemma to the set $X=\{x_1,...,x_k\}$ and the Whitney $(m-1)$-field $\VP^{(m,E)}[f]$. The Main Lemma \reff{X-SXN} produces a finite strictly increasing sequence $V=\{v_j\}_{j=1}^\ell$ in $E$ and a mapping $H:X\to 2^V$ which to every $x\in X$ assigns $m$ consecutive points
of $V$ possessing properties ($\blbig 1$)-($\blbig 5$) of the Main Lemma.
\par Using these objects, the sequence $V$ and the mapping $H$, in the next two lemmas we prove the required inequality \rf{VP-LM}.
\begin{lemma}\lbl{A2-LM} Let 
\bel{APL}
A^+=\smed_{j=1}^{k-1}\,\,
\smed_{i=0}^{m-1}\,\,
\frac{|L_{H(x_{j+1})}^{(i)}[f](x_{j})
-L_{H(x_{j})}^{(i)}[f](x_{j})|^p}
{(x_{j+1}-x_{j})^{(m-i)p-1}}.
\ee
\par Then $A^+\le C(m)^p\,\lambda^p$. (We recall that Assumption \reff{A-1} holds for the function $f$.)
\end{lemma}
\par {\it Proof.} Let $I_j$ be the smallest closed interval containing $H(x_j)\cup H(x_{j+1})$, $j=1,...,k-1$. Clearly, $\diam I_j=\diam (H(x_j)\cup H(x_{j+1}))$ and $x_{j+1},x_j\in I_j$ (because $x_j\in H(x_j)$ and $x_{j+1}\in H(x_{j+1})$, see property ($\blbig 2$) of the Main Lemma \reff{X-SXN}). Furthermore,
property ($\blbig 3$) of the Main Lemma \reff{X-SXN} implies that
\bel{LE-J}
I_j=[\min H(x_j),\max H(x_{j+1})].
\ee
\par Let us compare $\diam I_j$ with $x_{j+1}-x_{j}$ whenever $H(x_j)\ne H(x_{j+1})$. In this case Properties ($\blbig 2$) and ($\blbig 4$) of the Main Lemma \reff{X-SXN} tell us that $x_j\in H(x_j)$, $x_{j+1}\in
H(x_{j+1})$, and
$$
\diam H(x_j)+
\diam H(x_{j+1})\le 2(m+1)(x_{j+1}-x_j)\,.
$$
Hence,
\be
\diam I_j&=&\diam (H(x_j)\cup H(x_{j+1}))\le\diam H(x_j)+
\diam H(x_{j+1})+(x_{j+1}-x_j)\nn\\
&\le& 2(m+1)(x_{j+1}-x_j)+(x_{j+1}-x_j)=(2m+3)(x_{j+1}-x_j).
\nn
\ee
\par This inequality and definition \rf{APL} imply that
\bel{A-2}
A^+\le C(m)^p\,\sbig_{j=1}^{k-1}\,\,
\sbig_{i=0}^{m-1}\,\,
\frac{\max\limits_{I_j}|L_{H(x_{j+1})}^{(i)}[f]
-L_{H(x_{j})}^{(i)}[f]|^p}
{(\diam I_j)^{(m-i)p-1}}.
\ee
\par Note that each summand in the right hand side of  inequality \rf{A-2} equals zero if $H(x_j)=H(x_{j+1})$. Therefore, in our estimates of $A^+$, without loss of generality, we may assume that
\bel{H-IS}
H(x_j)\ne H(x_{j+1})~~~\text{for all}~~~j=1,...,k-1.
\ee
\par Property ($\blbig 1$) of the Main Lemma \reff{X-SXN} tells us that for every $j=1,...,k$, the set $H(x_j)$ consists of $m$ consecutive points of the strictly increasing sequence
\bel{V-DF2}
V=\{v_1,...,v_\ell\}\subset E.
\ee
\par Let $n_j$ be the index of the minimal point of $H(x_j)$ in the sequence $V$. Thus
$$
H(x_j)=\{v_{n_j},...,v_{n_j+m-1}\},~~~j=1,...,k-1.
$$
In particular, thanks to \rf{LE-J},
\bel{IJ-V}
I_j=[v_{n_j},v_{n_{j+1}+m-1}],~~~~j=1,..,k-1.
\ee
\par Let us apply Lemma \reff{LP-TSQ} (with $k=m-1$) to the sequence $\Yc=\{v_i\}_{i=n_j}^{n_{j+1}+m-1}$, the sets
$$
S_1=H(x_j)=\{v_{n_j},...,v_{n_j+m-1}\},~~~
S_2=H(x_{j+1})=\{v_{n_{j+1}},...,v_{n_{j+1}+m-1}\},
$$
and the closed interval $I=I_j=[v_{n_j},v_{n_{j+1}+m-1}]$. Also, let
$S^{(n)}_j=\{v_{n},...,v_{n+m}\},~~~n_j\le n\le n_{j+1}-1$.
\par In these settings Lemma \reff{LP-TSQ} tells us that
$$
\max\limits_{I_j}|L_{H(x_{j+1})}^{(i)}[f]
-L_{H(x_{j})}^{(i)}[f]|^p\le ((m+1)!)^p
\,(\diam I_j)^{(m-i)p-1}\,
\smed_{n=n_j}^{n_{j+1}}\,|\Delta^{m}f[S^{(n)}_j]|^p\,\diam S^{(n)}_j\,.
$$
This inequality and \rf{A-2} imply that
$$
A^+\le C(m)^p\,\smed_{j=1}^{k-1}\,
\smed_{n=n_j}^{n_{j+1}}\,(v_{n+m}-v_{n})\,
|\Delta^{m}f[v_{n},...,v_{n+m}]|^p.
$$
\par Let us prove that each interval $I_j$ contains at most $4m$ elements of the sequence $V$, i.e.,
\bel{IJ-VN}
n_{j+1}+m-n_j\le 4m~~~~ \text{for every}~~~j=1,..,k-1.
\ee
\par Indeed, let $x_j=v_{\vkp_j}$ and $x_{j+1}=v_{\vkp_{j+1}}$. Property ($\blbig 2$) of Main Lemma \reff{X-SXN} implies that
$0<\vkp_{j+1}-\vkp_{j}\le 2m$~ and~
$x_j\in H(x_j), x_{j+1}\in H(x_{j+1})$. But $\#H(x_{j})=\#H(x_{j+1})=m$ so that
$$
n_j\le\vkp_{j}\le n_j+m-1~~~\text{and}~~~ n_{j+1}\le\vkp_{j+1}\le n_{j+1}+m-1.
$$
\par These inequalities imply \rf{IJ-VN} because
$$
n_{j+1}+m-n_j\le n_{j+1}+m-\vkp_{j}+m-1\le
n_{j+1}-(\vkp_{j+1}-2m)+2m-1\le 4m.
$$
\par Property \rf{IJ-VN} enables us to estimate $A^+$ as follows: Let $\Ic=\{I_j:j=1,...,k-1\}$. Property ($\blbig 3$) of Main Lemma \reff{X-SXN} and \rf{H-IS} imply that
\bel{V-SI}
\{v_{n_j}\}_{j=1}^{k-1}~~~\text{is a {\it strictly increasing subsequence} of the sequence}~~V.
\ee
\par In turn, from \rf{IJ-V}, \rf{IJ-VN} and \rf{V-SI}, it follows that {\it  every interval $I_{j_0}\in\Ic$ has common points with at most $8m$ intervals $I_j\in \Ic$}. This leads us to the following property of the family $\Ic$: there exist subfamilies $\Ic_\nu\subset\Ic$, $\nu=1,...,\upsilon$, with $\upsilon\le 8m+1$, each consisting of {\it pairwise disjoint intervals}, such that
$\Ic=\cup\{\Ic_\nu:\nu=1,...,\upsilon\}$. The existence of subfamilies  $\{\Ic_\nu\}$ with these properties is immediate from the next well known statement from graph theory (see, e.g. \cite{JT}): {\it Every graph can be colored with one more color than the maximum vertex degree.}
\smsk
\par This, \rf{A-2} and inequality $\upsilon\le 8m+1$ imply that $A^+\le C(m)^p\,\max\{A_{\nu}:\nu=1,...,\upsilon\}$ where
$$
A_{\nu}=\,\smed_{j:I_j\in\Ic_\nu}\,\,
\smed_{n=n_j}^{n_{j+1}}\,(v_{n+m}-v_{n})\,
|\Delta^{m}f[v_{n},...,v_{n+m}]|^p.
$$
\par Since the intervals of each family $I_\nu$, $\nu=1,...,\upsilon$, are {\it pairwise disjoint}, the following inequality
$$
A_{\nu}\le \,\smed_{n=1}^{\ell-m}\,(v_{n+m}-v_{n})\,
|\Delta^{m}f[v_{n},...,v_{n+m}]|^p
$$
holds. (We recall that $\ell=\#V$, see \rf{V-DF2}.)
\par Applying Assumption \reff{A-1} to the right hand side of this inequality, we prove that $A_{\nu}\le \lambda^p$ for every $\nu=1,...,\upsilon$. Hence, $A^+\le C(m)^p\, \lambda^p$ completing the proof of the lemma.\bx
\begin{lemma}\lbl{FT-14} Inequality \rf{VP-LM} holds for every integer $k>1$ and every strictly increasing sequence $\{x_1,...,x_k\}\subset E$.
\end{lemma}
\par {\it Proof.} Let us replace the Hermite polynomials $\{P_{x_j}:j=1,...,k\}$ in the left hand side of inequality \rf{VP-LM} with the corresponding Lagrange polynomials $L_{H(x_j)}$. We will do this with the help of property ($\blbig 5$) of the Main Lemma \reff{X-SXN}.
\smsk
\par For every $j=1,...,k-1,$ and every $i=0,...,m-1,$ we have
\be
|P_{x_j}^{(i)}(x_j)-P_{x_{j+1}}^{(i)}(x_j)|&\le&
|P_{x_{j}}^{(i)}(x_{j})-L_{H(x_{j})}^{(i)}[f](x_{j})|+
|L_{H(x_{j})}^{(i)}[f](x_{j})
-L_{H(x_{j+1})}^{(i)}[f](x_{j})|
\nn\\
&+&
|L_{H(x_{j+1})}^{(i)}[f](x_{j})-P_{x_{j+1}}^{(i)}(x_{j})|.
\nn
\ee
Property ($\blbig 5$) of the Main Lemma (see \rf{FL-PH}) tells us that
$$
|P_{x_{j}}^{(i)}(x_{j})-L_{H(x_{j})}^{(i)}[f](x_{j})|+
|L_{H(x_{j+1})}^{(i)}[f](x_{j})-P_{x_{j+1}}^{(i)}(x_{j})|\le 2\ve
$$
proving that
$$
|P_{x_j}^{(i)}(x_j)-P_{x_{j+1}}^{(i)}(x_j)|\le
|L_{H(x_{j})}^{(i)}[f](x_{j})
-L_{H(x_{j+1})}^{(i)}[f](x_{j})|
+2\ve\,.
$$
Hence,
\bel{P-T1}
|P_{x_j}^{(i)}(x_j)-P_{x_{j+1}}^{(i)}(x_j)|^p\le
2^p\,|L_{H(x_{j})}^{(i)}[f](x_{j})
-L_{H(x_{j+1})}^{(i)}[f](x_{j})|^p+4^p\ve^p.
\ee
\par Let $A_1$ be the left hand side of inequality \rf{VP-LM}, i.e.,
\bel{A1-AW}
A_1=\smed_{j=1}^{k-1}\,\,
\smed_{i=0}^{m-1}\,\,
\frac{|P^{(i)}_{x_{j+1}} (x_j)-P^{(i)}_{x_j}(x_j)|^p}
{(x_{j+1}-x_{j})^{(m-i)p-1}},
\ee
and let
$$
A_2=4^p\,\smed_{j=1}^{k-1}\,\,\smed_{i=0}^{m-1}\,\,
(x_{j+1}-x_{j})^{1-(m-i)p}.
$$
\par Applying inequality \rf{P-T1} to each summand in the right hand side of \rf{A1-AW}, we obtain the following inequality:~ $A_1\le 2^p\,A^{+}+\ve^p\, A_2$. Recall that $A^+$ is the quantity defined by \rf{APL}. Lemma \reff{A2-LM} tells us that
$A^+\le C(m)^p\,\lambda^p$. Hence $A_1\le C(m)^p\,\lambda^p+\ve^p\, A_2$. But $\ve$ is an arbitrary positive constant, so that $A_1\le C(m)^p\,\lambda^p$ proving \rf{VP-LM} and the lemma.\bx
\bigskip

\par We are in a position to complete the proof of Theorem \reff{MAIN-TH}.
\medskip
\par {\it Proof of the sufficiency part of Theorem \reff{MAIN-TH}.} Definition \rf{VP-V} and Lemma \reff{FT-14} imply that
$$
\Nc_{m,p,E}\left(\VP^{(m,E)}[f]\right)\le C(m)\,\lambda=C(m)\,\NMP(f:E)\,.
$$
See \rf{NR-TR}. (We recall that $\lambda=\NMP(f:E)$.)  This inequality and the sufficiency part of Theorem \reff{JET-V} imply that 
$$
\|\VP^{(m,E)}[f]\|_{m,p,E}\le C(m)\, \Nc_{m,p,E}\left(\VP^{(m,E)}[f]\right)\le C(m)\,\NMP(f:E)\,.
$$
\par We recall that the quantity $\|\cdot\|_{m,p,E}$ is defined by \rf{N-VP}. This definition and the above inequality imply the existence of a function
$F\in\LMPR$ such that $T^{m-1}_x[F]=P_x$ on $E$, and
\bel{PEF-1}
\|F\|_{\LMPR}\le
2\,\|\VP^{(m,E)}[f]\|_{m,p,E}\le C(m)\,\NMP(f:E)\,.
\ee
\par We know that $P_x(x)=f(x)$ for every $x\in E$. (See \rf{FJ-H}-\rf{PX-M}.) Therefore,
$$
F(x)=T^{m-1}_x[F](x)=P_x(x)=f(x)~~~\text{for all}~~~x\in E\,.
$$
\par Thus, $F\in\LMPR$ and $F|_E=f$ proving that $f\in \LMPR|_E$. Furthermore, definition \rf{N-LMPR} and inequali\-ty \rf{PEF-1} imply that
$$
\|f\|_{\LMPR|_E}\le \|F\|_{\LMPR}\le C(m)\,\NMP(f:E)
$$
proving the sufficiency.\bx
\smsk
\par The proof of Theorem \reff{MAIN-TH} is complete.\bx
\msk
\begin{remark} {\em Given a function $f$ on $E$, let us indicate the main steps of our extension algorithm suggested in Sections 3 and 4:
\smsk
\par {\it Step 1.} We construct the family of sets
$\{\SH_x: x\in E\}$ and the family of points $\{s_x: x\in E\}$ satisfying conditions of Proposition \reff{SET-SX} and Proposition \reff{PR-SX};
\smsk
\par {\it Step 2.} At this step we construct
the Whitney $(m-1)$-field $\VP^{(m,E)}[f]=\{P_x\in\PMO:x\in E\}$ satisfying conditions (i), (ii) of Definition
\reff{P-X};
\smsk
\par {\it Step 3.} We define the extension $F$ by formula \rf{DEF-F}.
\msk
\par We denote the extension $F$ by $F=\EXT_E(f:\LMPR)$.
Clearly, $F$ depends on $f$ {\it linearly} proving that $\EXT_E(\cdot:\LMPR)$ is a {\it linear extension operator}. Theorem \reff{MAIN-TH} states that its operator norm is bounded by a constant depending only on $m$.}\rbx
\end{remark}

\begin{remark} {\em In \cite[Section 4.4]{Sh-LV-2018} we give an alternative proof of Theorem \reff{DEBOOR} based on the extension algorithm described in Section 3 and Section 4 of the present paper.
\par We note that, for the case of sequences, this extension method can be simplified considerably. Indeed, we prove in \cite[Remark 3.13]{Sh-LV-2018} that for each $x\in E$, the set $\SH_x$ consist of $m$ consecutive terms of the sequence $E$. In this case part (ii) of Definition \reff{P-X} tells us that $P_x$ coincides with the Lagrange polynomial $L_{\SH_x}$ interpolating $f$ on $\SH_x$.
\par  In turn, this property immediately implies a variant of the Main Lemma for sequences, see
\cite[Lemma 4.13]{Sh-LV-2018}, where we set $H(x)=\SH_x$, for every $x\in E$. The required properties of the sets $\{H(x):x\in E\}$ for this case are immediate from Proposition \reff{SET-SX}.
\par We also note that if $E=\{x_i\}_{i=\ell_1}^{\ell_2}$ is a strictly increasing sequence of points in $\R$, and $f$ is a function on $E$, the extension $F=\EXT_E(f:\LMPR)$ is a piecewise polynomial $C^{m-1}$-function which coincides with a polynomial of degree at most $2m-1$ on each interval $(x_i,x_{i+1})$. This enables us to reformulate this property of $F$ in terms of Spline Theory as follows: {\it The extension $F$ is an interpolating $C^{m-1}$-smooth spline of order\, $2m$\, with knots\, $\{x_i\}_{i=\ell_1}^{\ell_2}$.}
\par Details are spelled out in \cite{Sh-LV-2018}.
}\rbx
\end{remark}

\SECT{7. The Finiteness Principle for $L^m_\infty(\R)$ traces: multiplicative finiteness constants.}{7}
\addtocontents{toc}{7. The Finiteness Principle for $L^m_\infty(\R)$ traces: multiplicative finiteness constants. \hfill\thepage\par \VST}

\indent\par Let $m\in\N$. Everywhere in this section we assume that $E$ is a closed subset of $\R$ with $\#E\ge m+1$.
\par We will discuss equivalence \rf{T-R1} which states that
\bel{FP-EQ}
\|f\|_{\LMIR|_E}\sim
\sup_{S\subset E,\,\,\# S=m+1}
|\Delta^mf[S]|
\ee
for every function $f\in\LMIR|_E$. The constants in this equivalence depend only on $m$.
\par We can interpret this equivalence as a special case of the following {\it Finiteness Principle} for the space $\LMIR$.
\begin{theorem}\lbl{FP-LM} Let $m\in\N$. There exists a constant $\gamma=\gamma(m)>0$ depending only on $m$, such that the following holds: Let $E\subset\R$ be a closed set, and let $f:E\to\R$.
\par For every subset $E'\subset E$ with at most $N=m+1$ points, suppose there exists a function $F_{E'}\in\LMIR$ with the seminorm $\|F_{E'}\|_{\LMIR}\le 1$, such that $F_{E'}=f$ on $E'$.
\par Then there exists a function $F\in\LMIR$ with the seminorm $\|F\|_{\LMIR}\le \gamma$ such that $F=f$ on $E$.
\end{theorem}
\par {\it Proof.} The result is immediate from equivalence \rf{FP-EQ} and definition \rf{N-LMPR}. Details are spelled out in \cite[Section 5]{Sh-LV-2018}.\bx
\msk
\par We refer to the number $N=m+1$ as {\it a finiteness number} for the space $\LMIR$. Clearly, the value $N(m)=m+1$ in the finiteness Theorem \reff{FP-LM} is sharp; in other words, Theorem \reff{FP-LM} is false in general if $N=m+1$ is replaced by some number $N<m+1$.
\par We note that the Finiteness Principle also holds for the space $\LMRN$ for all $m,n\in\N$; in this case a corresponding number $N$ and a constant $\gamma$ depend only on $n$ and $m$. See \cite{Sh-1987} for the case $m=2$, $n\in\N$, and \cite{F2} for the general case of $m,n\in\N$. It is also shown in \cite{F2} that a similar finiteness principle holds for the space $\WMRN$.
\smallskip
\par Theorem \reff{FP-LM} implies the following: For every function $f\in\LMIR|_E$, we have
\bel{MFC-G}
\|f\|_{\LMIR|_E}\le \gamma(m)
\,\sup_{S\subset E,\,\,\# S= m+1}\,\|f|_S\|_{\LMIR|_S}.
\ee
(Clearly, the converse inequality is trivial and holds with $\gamma(m)=1$.) We refer to any constant $\gamma=\gamma(m)$ which satisfies \rf{MFC-G} as {\it a multiplicative finiteness constant} for the space $\LMIR$.
\par The following natural question arises:
\begin{question}\lbl{Q-1} What is the sharp value of the multiplicative finiteness constant for $\LMIR$, i.e., the infimum of all multiplicative finiteness constants for $\LMIR$ for the finiteness number $N=m+1$?
\end{question}
\par We denote this sharp value of $\gamma(m)$ by $\gmr$. Thus,
\bel{D-GM}
\gmr=\sup\,\,\frac{\|f\|_{\LMIR|_E}}
{\sup\{\|f|_S\|_{\LMIR|_S}: S\subset E,\,\# S= m+1\}}
\ee
where the supremum is taken over all closed sets $E\subset\R$ with $\#E\ge m+1$, and all $f\in\LMIR|_E$.
\smallskip
\par The next theorem answers to Question \reff{Q-1} for $m=1,2$ and provides lower and upper bounds for $\gmr$ for $m>2$. These estimates show that $\gmr$ {\it grows exponentially} as $m\to +\infty$.
\begin{theorem}\lbl{G-SH} (i) $\gsho=1$ and $\gsht=2$;
\smallskip
\par (ii) For every $m\in\N$, $m>2$, the following inequalities
$$
\left(\frac{\pi}{2}\right)^{m-1}<\gmr<(m-1)\,9^m
$$
hold.
\end{theorem}
\par The proof of this theorem relies on works \cite{Fav,deB4,deB5} devoted to calculation of a certain constant $K(m)$ related to optimal extensions of $\LMIR$-functions. This constant is defined by
\bel{D-KM}
K(m)=\sup\,\,\frac{\|f\|_{\LMIR|_X}}
{\max\{m!\,|\Delta^mf[x_i,...,x_{i+m}]|: i=1,...,n\}}
\ee
where the supremum is taken over all $n\in\N$, all  strictly increasing sequences $X=\{x_1,...,x_{m+n}\}\subset\R$ and all functions $f\in\LMIR|_X$.
\par The constant $K(m)$ was introduced by Favard \cite{Fav}. (See also \cite{deB4,deB5}.) Favard \cite{Fav} proved that $K(2)=2$, and de Boor found efficient lower and upper bounds for $K(m)$.
\smsk
\par We prove that $\gmr=K(m)$, see Lemma \reff{GM-KM} below. This formula and aforementioned results of Favard and de Boor imply the required lower and upper bounds for $\gmr$ in Theorem \reff{G-SH}.
\smallskip
\par We will need a series of auxiliary lemmas.
\begin{lemma}\lbl{DD-M} Let $S\subset\R$, $\#S=m+1$, and let $f:S\to\R$. Then
$\|f\|_{\LMIR|_S}=m!\,|\Delta^mf[S]|$.
\end{lemma}
\par {\it Proof.} Let $A=\min S$, $B=\max S$. Let $F\in\LMIR$ be an arbitrary function such that $F|_S=f$.
Inequality \rf{F-LMIR} tells us that
$m!\,|\Delta^mf[S]|=m!\,|\Delta^mF[S]|\le\|F\|_{\LMIR}$.
\par Taking the infimum in this inequality over all functions $F\in\LMIR$ such that $F|_S=f$, we obtain the inequality $m!\,|\Delta^mf[S]|\le \|f\|_{\LMIR|_S}$.
\par Let us prove the converse inequality. Let $F=L_S[f]$. Then
$$
m!\,|\Delta^{m}f[S]|=|L^{(m)}_S[f]|=\|F\|_{\LMIR}~~~~~
\text{(see \rf{D-LAG})}.
$$
But $F|_S=L_S[f]|_S=f$, so that
$
\|f\|_{\LMIR|_S}\le \|F\|_{\LMIR}=m!\,|\Delta^mf[S]|
$
proving the lemma.\bx
\begin{lemma}\lbl{CR-1} Let
$E=\{x_1,...,x_{m+n}\}\subset\R$, $n\in\N$, be a strictly increasing sequence, and let $f$ be a function on $E$. Then
$$
\max_{S\subset E,\,\#S=m+1}|\Delta^mf[S]|=
\max_{i=1,...,n}\,|\Delta^mf[x_i,...,x_{i+m}]|\,.
$$
\end{lemma}
\par {\it Proof.} The lemma is immediate from the following property of divided differences (see
\cite[p. 15]{FK-88}): Let $S\subset E$, $\#S=m+1$. There exist $\alpha_i\ge 0$, $i=1,...,n$ with $\alpha_1+...+\alpha_{n}=1$ such that
$$
\Delta^mf[S]=\smed_{i=1}^{n}\,\,\alpha_i\,
\Delta^mf[x_i,...,x_{i+m}]\,.\BX
$$
\begin{lemma}\lbl{R-FS} Let $E\subset\R$ be a closed set,
and let $f\in\LMIR|_E$. Then
$$
\|f\|_{\LMIR|_E}=\sup_{E'\subset E,\,\#E'<\infty}\|f|_{E'}\|_{\LMIR|_{E'}}\,.
$$
\end{lemma}
\par {\it Proof.} We recall the following well known fact: for every closed bounded interval $I\subset\R$ a ball in the space  $L^m_\infty(I)$ is a precompact subset in the space $C(I)$. The lemma readily follows from this statement. We leave the details to the interested reader.\bx
\begin{lemma}\lbl{GM-KM} For every $m\in\N$ the following equality $\gmr=K(m)$ holds.
\end{lemma}
\par {\it Proof.} The inequality $K(m)\le\gmr$ is immediate from Lemma \reff{DD-M}, Lemma \reff{CR-1} and definition \rf{D-KM}. The converse inequality directly follows from Lemma \reff{R-FS} and definition \rf{D-GM}.
\par Details are spelled out in \cite[Section 6] {Sh-LV-2018}.\bx
\medskip
\par {\it Proof of Theorem \reff{G-SH}.} The equality
$\gsho=1$ is immediate from the well known fact that a function satisfying a Lipschitz condition on a subset of $\R$ can be extended to all of $\R$ with preservation of the Lipschitz constant.
\par As we have mentioned above, $K(2)=2$ (Favard  \cite{Fav}). de Boor \cite{deB4,deB5} proved that
\bel{DB-R12}
\left(\frac{\pi}{2}\right)^{m-1}<\cm\le K(m)\le \CMM<(m-1)\,9^m~~~\text{for each}~~~m>2.
\ee
Here
$$
\cm= \left(\frac{\pi}{2}\right)^{m+1}
\left/\,\,
\left
(\,\smed_{j=-\infty}^\infty\,\,((-1)^j/(2j+1))^{m+1}
\right)
\right.~~~~\text{and}~~~~
\CMM= (2^{m-2}/m)+\,\smed_{i=1}^m\, \binom{m}{i}\,\binom{m-1}{i-1}\,4^{m-i}\,.
$$
\par On the other hand, Lemma \reff{GM-KM} tells us that $\gmr=K(m)$ which together with the equality $K(2)=2$ and inequalities \rf{DB-R12} implies statements (i) and (ii) of the theorem.
\par The proof of Theorem \reff{G-SH} is complete.\bx
\msk
\par We finish the section with two remarks.
\begin{remark} {\em One can readily see that $\theta_1=1$ and $\theta_2=2$. This together with part (i) of Theorem \reff{G-SH} implies that $\cm=K(m)$ for $m=1,2$. 
\par In turn, thanks to \rf{DB-R12}, $\cm\le K(m)$ for $m>2$. de Boor \cite{deB4} proved this inequality by showing that for a function $f(i)=(-1)^i$ defined on $\Z$, the following two properties hold:
\par (i) $m!\,|\Delta^mf[i,...,i+m]|=2^m$ for every $i\in\Z$; ~(ii)  $\|f\|_{\LMIR|_{\Z}}=\cm\,2^m$.
\smsk
\par Furthermore, it is shown that {\it the Euler spline}
$\Ec_m:\R\to\R$ (a classical object of Spline Theory introduced by Schoenberg \cite{Sch4}) provides an optimal (in $\LMIR$) extension of $f$ from $\Z$ to all of $\R$.
\par We refer the reader to \cite[Section 3]{deB4}, \cite[Lecture 4]{Sch4}, \cite[Section 6.2]{Sh-LV-2018} for more detail.\rbx}
\end{remark}

\begin{remark} {\em  In fact we know very little about the values of the constants $\gamma^\sharp(\cdot)$ and its analogs for the spaces $\LMRN$ and $\WMRN$.
Apart from the results related to $\gmr$ which we present in this section, there is perhaps only one other result in this direction. It is due to Fefferman and Klartag \cite{FK}. They studied the behavior of the sharp multiplicative finiteness constant for the space $\WTRT$ defined by
$$
\gamma^\sharp(N;\WTRT)=\sup\,\,\frac{\|f\|_{\WTRT|_E}}
{\sup\{\|f|_S\|_{\WTRT|_S}: S\subset E,\,\# S\le N\}}.
$$
\par It is proven in \cite{FK} that, after a certain renormalization of the space $\WTRT$, the following statement holds: {\it there exists an absolute constant $c>0$ such that $\gamma^\sharp\left(N;\WTRT\right)>1+c$
for every positive integer $N$.}\rbx}
\end{remark}



\begin{thebibliography}{ABCD}

\addtocontents{toc}{References \hfill \thepage\par}
\bibitem {BZ} I. S. Berezin, N. P. Zhidkov, Computing methods, v. 1 and 2, Pergamon Press, London, 1965.
\bibitem {B} Yu. Brudnyi, Criteria for the existence of derivatives in $L_p$. (Russian) Mat. Sb. (N.S.),  73 No. 115  (1967) 42--64.
\bibitem {B2} Yu. A. Brudnyi, Spaces that are definable by means of local approximations. (Russian) Trudy Moscov. Math. Obshch. 24 (1971) 69--132; English transl.: Trans. Moscow Math. Soc. 24 (1974) 73--139.
\bibitem {BS1} Yu. Brudnyi, P. Shvartsman,
Generalizations  of Whitney  Extension Theorem,
Intern. Math.  Research Notices, No. 3, (1994) 129--139.
\bibitem {BS2} Yu. Brudnyi, P. Shvartsman,
The Whitney Problem  of  Existence of a Linear
Extension Operator, J. Geom. Anal., {7}, No. 4, (1997) 515--574.
\bibitem {BS3} Yu. Brudnyi, P. Shvartsman,
Whitney Extension Problem for Multivariate
$C^{1,\omega}$-functions, Trans. Amer. Math. Soc. {353} No. 6,  (2001) 2487--2512.
\bibitem {C1} A. P. Calder\'{o}n, Estimates for
    singular integral operators in terms of maximal
    functions, Studia Math. 44 (1972) 563--582.
\bibitem {CS} A. P. Calder\'{o}n, R. Scott, Sobolev
    type inequalities for $p>0$,  Studia Math. 62 (1978)
    75--92.
\bibitem {ChS} C. K. Chui and P. W. Smith, On $H^{m,\infty}$-splines, SIAM J. Numer. Anal. 11 (1974) 554--558.
\bibitem {CSW1} C. K. Chui, P. W. Smith, and J. D. Ward, Favard's solution is the limit of $W^{k,p}$-splines, Trans. Amer. Math. Soc. 220 (1976) 299--305.
\bibitem {deB4} C. de Boor, How small can one make the derivatives of an interpolating function? J. Approx.
Theory, 13 (1975) 105--116.
\bibitem {deB5} C. de Boor, A smooth and local interpolant with $k$-th derivative, In: Numerical Solutions of Boundary Value Problems for Ordinary Differential Equations, (Proc. Sympos., Univ. Maryland, Baltimore, Md., 1974),  pp. 177--197. Academic Press, New York, 1975.
\bibitem {deB2} C. de Boor, On "best" interpolation, J. Approx. Theory 16 (1976) 28--42.
\bibitem {deB6} C. de Boor, Splines as linear combinations of B-splines. A survey,  Approximation theory, II (Proc. Internat. Sympos., Univ. Texas, Austin, Tex., 1976),  pp. 1--47. Academic Press, New York, 1976.
\bibitem {DL} R. DeVore, G. Lorentz, Constructive approximation, Fundamental Principles of Mathema\-ti\-cal Sciences, 303. Springer-Verlag, Berlin, 1993.
\bibitem {Es} D. Est\'evez, Explicit traces of functions from Sobolev spaces and quasi-optimal linear interpolators, Math. Inequal. Appl. 20 (2017) 441--457.
\bibitem {Fav} J. Favard, Sur l'interpolation, J. Math. Pures Appl. 19 (1940) 281--306.
\bibitem {F2} C. Fefferman, A sharp form of Whitney extension theorem, Annals of Math. 161, No. 1, (2005) 509--577.
\bibitem {F3} C. Fefferman,  $C^m$ Extension by Linear Operators, Annals of Math. 166, No. 3, (2007) 779--835.
\bibitem {F-IO} C. Fefferman, Extension of $C^{m,\omega}$-Smooth Functions by Linear Operators, Rev. Mat. Iberoamericana 25, No. 1, (2009) 1--48.
\bibitem {F-Bl} C. Fefferman, Whitney extension problems and interpolation of data, Bulletin A.M.S. 46, No. 2, (2009) 207--220.
\bibitem {F9} C. Fefferman, The $C^m$ norm of a function with prescribed jets I., Rev. Mat. Iberoamericana 26 (2010) 1075--1098.
\bibitem {FIL} C. Fefferman, A. Israel, G. K. Luli,
Sobolev extension by linear operators, J. Amer. Math. Soc. 27 (2014) 69--145.
\bibitem {FK} C. Fefferman, B. Klartag, An example related to Whitney extension with almost minimal $C^m$ norm,  Rev. Mat. Iberoamericana 25 (2009) 423--446.
\bibitem {FJ1} S. Fisher, J. Jerome, Minimum Norm Extremals in Function Spaces, Lect. Notes Math., 479 (1975).
\bibitem {FK-88} A. Fr\"olicher, A. Kriegl, Linear spaces and differentiation theory, Pure and Applied Mathematics, J. Wiley, Chichester, 1988.
\bibitem {FK-93} A. Fr\"olicher, A. Kriegl, Differentiable extensions of functions, Differ. Geom. Appl. 3 (1993)
71--90.
\bibitem {Gol} M. Golomb, $H^{m,p}$-extensions by $H^{m,p}$-splines, J. Approx. Theory 5 (1972) 238--275.
\bibitem {Is} A. Israel, A Bounded Linear Extension Operator for $L^{2,p}(\R^2)$, Annals of Math.  178 (2013) 1--48.
\bibitem {JR} A. Jakimovski, D. C. Russell, On an interpolation problem and spline functions.  General inequalities, 2 (Proc. Second Internat. Conf., Oberwolfach, 1978),  pp. 205--231, Birkhäuser, Basel-Boston, Mass., 1980.
\bibitem {JT} T. R. Jensen, B. Toft, Graph coloring problems, Wiley-Interscience Series in Discrete Mathematics and Optimization. A Wiley-Interscience Publication. John Wiley $\&$ Sons, Inc., New York, 1995.
\bibitem {JS} J. W. Jerome, L. L. Schumaker,
Characterizations of Functions with Higher Order Derivatives in $L_p$, Trans. Amer. Math. Soc. 143 (1969) 363--371.
\bibitem {Jon} A. Jonsson, The Trace of the Zygmund Class $\Lambda(\R)$ to Closed Sets and Interpolating Polynomials, J. Approx. Th. 4 (1985) 1--13.
\bibitem {KM} A. Kriegl, P. W. Michor, The convenient setting of global analysis, Math. surveys and monographs, v. 53, AMS, 1997.
\bibitem {Kun2}  T. Kunkle, Favard's interpolation problem in one or more variables, Constr. Approx.  18  (2002),  no. 4, 467--478.
\bibitem {L} G.K. Luli, Whitney Extension for $W^k_p(E)$ in One Dimension, (2008). Notes available at
https://www.math.ucdavis.edu/~kluli/LmpNotes.pdf
\bibitem {Mer} J. Merrien, Prolongateurs de fonctions differentiables d'une variable r\'eelle, J. Math. Pures et Appl.  45 (1966) 291--309.
\bibitem {P2} A. Pinkus, Uniqueness of Smoothest Interpolants, East J. Approx. 3 (1997) 377--380.
\bibitem {R} F. Riesz, Untersuchungen \"uber  Systeme integrierbarer Funktionen, Math. Ann. 69 (1910) 449--497.
\bibitem {Sch2} I. J. Schoenberg, Spline interpolation and the higher derivatives, Proc. Nar. Acad. Sci. U.S.A. 51 (1964) 24--28.
\bibitem {Sch4} I. J. Schoenberg, Cardinal Spline Interpolation, CBMS, SIAM, Philadelphia, 1973.
\bibitem {SHE1} I. A. Shevchuk, Extension of functions which are traces of functions belonging to $H^{\varphi}_k$ on arbitrary subset of the line, Anal. Math., No. 3, (1984) 249--273.
\bibitem {Sh-1987} P. Shvartsman, On the traces of functions of the Zygmund class, Sib. Mat. Zh. 28 (1987) 203--215; English transl. in Sib. Math. J. 28 (1987) 853--863.
\bibitem {Sh-2008} P. Shvartsman, The Whitney extension problem and Lipschitz selections of set-valued mappings in jet-spaces, Trans. Amer. Math. Soc.  360 (2008)  5529--5550.
\bibitem {Sh2} P. Shvartsman, Sobolev $W^1_p$-spaces on closed subsets of $\RN$, Advances in Math. 220 (2009) 1842--1922.
\bibitem {Sh3} P. Shvartsman, Sobolev $L^2_p$-functions on closed subsets of $\R^2$, Advances in Math. 252 (2014) 22--113.
\bibitem {Sh5} P. Shvartsman, Whitney-type extension theorems for jets generated by Sobolev functions,  Advances in Math. 313 (2017) 379--469.
\bibitem {Sh-LV-2018} P. Shvartsman, Sobolev functions on closed subsets of the real line: long version, arXiv:1808.01467
\bibitem {Sh-2018} P. Shvartsman, Sobolev functions on closed subsets of the real line, arXiv:1710.07826
\bibitem {Sm-2} P. W. Smith, $W^{r,p}(R)$-splines, J. Approx. Th. 10 (1974) 337--357.
\bibitem {W1}  H. Whitney, Analytic extension of differentiable functions defined in closed sets, Trans.  Amer. Math. Soc.  36 (1934) 63--89.
\bibitem {W2} H. Whitney, Differentiable functions defined in closed sets. I., Trans. Amer. Math. Soc. 36 (1934) 369--387.
\end{thebibliography}
\end{document}